\documentclass[a4paper,12pt,reqno]{amsart}

\usepackage{amssymb,amsmath,amsthm,latexsym}
\usepackage{amscd}
\makeatletter

\@addtoreset{equation}{section}

\@addtoreset{figure}{section}
\makeatother

\theoremstyle{plain}
\newtheorem{Theorem}{Theorem}[section]

\newtheorem{Proposition}[Theorem]{Proposition}
\newtheorem{Proposition-Definition}[Theorem]{Proposition-Definition}
\newtheorem{Lemma}[Theorem]{Lemma}

\theoremstyle{definition} 
\newtheorem{Definition}[Theorem]{Definition}
\newtheorem{Remark}[Theorem]{Remark}
\newtheorem{Example}[Theorem]{Example}

\newtheorem{Observation}[Theorem]{Observation}

\newenvironment{acknowledgments}[1][Acknowledgments] {\\\textit{#1.~~}\it}{}

\def\sO{{\mathcal{O}}}
\def\J{{\mathcal{J}}} 
\def\Z{{\mathbb{Z}}}

\def\Q{{\mathbb{Q}}} 
\def\R{{\mathbb{R}}} 
\def\C{{\mathbb{C}}}

\def\x{{\boldsymbol{x}}}
\def\t{{\boldsymbol{t}}}
\def\f{{\boldsymbol{f}}}
\def\vv{{\boldsymbol{v}}}
\def\w{{\boldsymbol{w}}}
\def\a{{\mathfrak{a}}}

\def\m{{\mathfrak{m}}}
\def\p{{\mathfrak{p}}}
\def\q{{\mathfrak{q}}}
\def\Spec{\mathop{\mathrm{Spec}}\nolimits}
\def\Ann{\mathop{\mathrm{Ann}}\nolimits}

\def\max{\mathop{\mathrm{max}}\nolimits}

\def\lct{\mathop{\mathrm{lct}}\nolimits}

\def\gr{\mathop{\mathrm{Gr}}\nolimits}

\def\initial{\mathop{\mathrm{in}}\nolimits}
\def\ord{\mathop{\mathrm{ord}}\nolimits}
\begin{document}
\title{Algorithms for computing multiplier ideals}
\author{Takafumi Shibuta
}
\address{
\begin{flushleft}
	\hspace{0.3cm} Department of Mathematics, Rikkyo University, Nishi-Ikebukuro, Tokyo 171-8501, Japan\\
	\hspace{0.3cm} JST, CREST, Sanbancho, Chiyoda-ku, Tokyo, 102-0075, Japan\\
\end{flushleft}
}
\email{shibuta@rikkyo.ac.jp}
\date{}
\baselineskip 15pt
\footskip = 32pt
\begin{abstract}
We give algorithms for computing multiplier ideals using Gr\"obner bases in Weyl algebras. 
The algorithms are based on a newly introduced notion which is a variant of Budur--Musta\c{t}\v{a}--Saito's (generalized) Bernstein--Sato polynomial. 
We present several examples computed by our algorithms. 
\end{abstract}
\maketitle
\tableofcontents
\section{Introduction}
Multiplier ideals are an important tool in higher-dimensional algebraic geometry, and one can view them as measuring singularities. 
They are defined as follows:
let $X$ be a smooth complex variety and $\a \subseteq \sO _X$ be an ideal sheaf of $X$. 
Suppose that $\pi:\widetilde{X} \to X$ is a log resolution of $\a$, that is, $\pi$ is a proper birational morphism, 
$\widetilde{X}$ is smooth and $\pi^{-1}V(\a)=F$ is a divisor with simple normal crossing support. 
If $K_{\widetilde{X}/X}$ is the relative canonical divisor of $\pi$, then the multiplier ideal of $\a$ with exponent $c \in \R_{\ge 0}$ is 
\[
\J(\a^c)=\J(c \cdot \a)=\pi_*\sO_{\widetilde{X}}(K_{\widetilde{X}/X}-\lfloor cF \rfloor) \subseteq \sO_X.
\]
A positive rational number $c$ is called a {\it jumping coefficient} if $\J(\a^{c-\varepsilon})\neq \J(\a^{c})$ for all $\varepsilon>0$, 
and the minimal jumping coefficient is called the {\it log-canonical threshold} of $\a$ and denoted by $\lct(\a)$. 
Since multiplier ideals are defined via log resolutions, it is difficult to compute them in general 
(when the ideal $\a$ is a monomial ideal or a principal ideal generated by a non-degenerate polynomial, 
there is a combinatorial description of the multiplier ideals $\J(\a^c)$. See \cite{H} and \cite{H2}). 
In this paper, we will give algorithms for computing multiplier ideals using the theory of $D$-modules. 


The Bernstein-Sato polynomial (or $b$-function) is one of the main objects in the theory of $D$-modules. 
It has turned out that jumping coefficients are deeply related to Bernstein-Sato polynomials.
For a given polynomial $f \in \C[x_1, \dots, x_n]$, 
the Bernstein-Sato polynomial $b_f(s)$ of $f$ is the monic polynomial in one variable $b(s) \in \C[s]$ of minimal degree 
having the property that there exists a linear differential operator $P(x,s)$ such that $b(s) f^s =P(x,s)f^{s+1}$. 
Koll\'ar \cite{Ko} proved that the log canonical threshold of $f$ is the minimal root of $b_f(-s)$. 
Furthermore, Ein--Lazarsfeld--Smith--Varolin \cite{ELSV} extended Koll\'ar's result to higher jumping coefficients: 
they proved that all jumping coefficients in the interval $(0,1]$ are roots of $b_f(-s)$. 
Recently Budur--Musta\c{t}\v{a}--Saito introduced the notion of Bernstein-Sato polynomials of arbitrary ideal sheaves 
using the theory of $V$-filtrations of Kashiwara \cite{Ka} and Malgrange \cite{M}. 
They then gave a criterion for membership of multiplier ideals in terms of their Bernstein-Sato polynomials, 
and proved that all jumping coefficients of $\a\subset \sO_X$ in the interval $(\lct(\a),\lct(\a)+1]$ 
are roots of the Bernstein-Sato polynomial of $\a$ up to sign. 

It is difficult to compute Bernstein-Sato polynomials in general, 
but Oaku \cite{O1}, \cite{O2}, \cite{O3} gave algorithms for computing Bernstein-Sato polynomials $b_f(s)$ 
using Gr\"obner bases in Weyl algebras 
(algorithms for computing Gr\"obner bases in Weyl algebras are implemented in some computer systems, 
such as Kan/Sm1 \cite{T} and Risa/Asir \cite{N}). 
In this paper, we give algorithms for computing Budur--Musta\c{t}\v{a}--Saito's Bernstein-Sato polynomials 
(Theorems \ref{Algorithm for global b-functions 1}, \ref{Algorithm for global b-functions 2} and \ref{Algorithm for local b-functions}). 
Our algorithms are natural generalizations of Oaku's algorithm. 

The other ingredient of this paper is algorithms for computing multiplier ideals. 
The algorithm for computing generalized Bernstein-Sato polynomials enables us to solve the membership problem for multiplier ideals, 
but does not give a system of generators of multiplier ideals. 
We modify the definition of Budur--Musta\c{t}\v{a}--Saito's Bernstein-Sato polynomials 
to determine a system of generators of the multiplier ideals of a given ideal (Definition \ref{new Bernstein-Sato polynomial}). 
Then we obtain algorithms for computing our Bernstein-Sato polynomials 
and algorithms for computing multiplier ideals (Theorem \ref{algorithm for computing multiplier ideals 1}, Theorem \ref{algorithm for computing multiplier ideals 2}). 
Our algorithms are based on the theory of Gr\"obner bases in Weyl algebras 
(see \cite{OT} and \cite{SST} for a review of Gr\"obner bases in Weyl algebras and their applications). 
We conclude the paper by presenting several examples computed by our algorithms. \\
\begin{acknowledgments}\upshape
The author thanks Shunsuke Takagi for useful comments and discussions. 
He also would like to thank Naoyuki Shinohara and Kazuhiro Yokoyama for warm encouragement and support. 
\end{acknowledgments}
\section{Preliminaries}
\subsection{Gr\"obner bases in Weyl algebras}
We denote by $\C$ the complex number field. 
When we use a computer algebra system, we may work with a computable field $\Q(z_1,\dots,z_l)\subset \C$ 
which is sufficient to express objects that appear in the computations.  
Let $X$ be the affine space $\C^{d}$ with the coordinate system $\x= (x_1,\dots,x_d)$, 
and $\C[\x]=\C[x_1,\dots,x_d]$ a polynomial ring over $\C$ which is the coordinate ring of $X$. 
We denote by ${\partial_\x} = (\partial_{x_1} ,\dots, \partial_{x_d})$ the partial differential operators where 
$\partial_{x_i} = \frac{\partial}{\partial{x_i}}$. 
We set 
\[
D_X=\C\langle \x,{\partial_\x}\rangle=\C\langle x_1,\dots,x_d,\partial_{x_1},\dots,\partial_{x_d}\rangle, 
\]
the rings of differential operators of $X$, and call it the {\it Weyl algebra} (in $d$ variables). 
This ring is non-commutative $\C$-algebra with the commutation rules 
\[
x_ix_j=x_jx_i,~~\partial_{x_i}\partial_{x_j}=\partial_{x_j}\partial_{x_i},
~~\partial_{x_i}x_j=x_j\partial_{x_i}~~\mathrm{for}~~i\neq j,~~\mathrm{and}~\partial_{x_i}x_i=x_i\partial_{x_i}+1.
\]
We write $\langle P_1,\dots,P_r\rangle$ for the left ideal of $D_X$ generated by $P_1,\dots,P_r\in D_X$. 
We use the notation 
$\x^{\boldsymbol{\mu}}=\prod_{i=1}^d x_i^{{\mu}_{i}}$,  
and 
${\partial_\x}^{\boldsymbol{\nu}}=\prod_{i=1}^d\partial_{x_i}^{{\nu}_{i}}$ 
for $\boldsymbol{\mu} = ({\mu}_{1}, \dots , {\mu}_{d})$, $\boldsymbol{\nu}=({\nu}_{1},\dots,{\nu}_{d})\in\Z_{\ge 0}^d$. 
We denote by $|\boldsymbol{\mu}| := {\mu}_{1} + \dots+ {\mu}_{d}$ the side of $\mu$. 
We call a real vector $(\vv,\w)=(v_{1},\dots,v_{d},w_{1},\dots,w_{d})\in \R^{d}\times\R^{d}$ a {\it weight vector} if 
\[
v_i+w_i\ge 0 \mathrm{~for~}i=1,2,\dots,d.
\]
We define the ascending filtration $\cdots \subset F^{(\vv,\w)}_0{D_X}\subset F^{(\vv,\w)}_1{D_X}\subset \cdots$ on $D_X$ 
with respect to the weight vector $(\vv,\w)$ by 
\[
F^{(\vv,\w)}_m{D_X}=\{\sum_{\vv\cdot \boldsymbol{\mu}+\w\cdot \boldsymbol{\nu}\le m} 
a_{\boldsymbol{\mu\nu}}\x^{\boldsymbol{\mu}}{\partial_\x}^{\boldsymbol{\nu}}\mid 
a_{\boldsymbol{\mu\nu}}\in \C\}\subset D_X 
\]
where $\vv\cdot \boldsymbol{\mu}=\sum v_i\mu_i$ is the usual inner product of $\vv$ and $\boldsymbol{\mu}$. 
Then we have 
\[
F^{(\vv,\w)}_{m_1}{D_X}\cdot F^{(\vv,\w)}_{m_2}{D_X}\subset F^{(\vv,\w)}_{m_1+m_2}{D_X}
\] 
for all $m_1$, $m_2\in \Z$ by the conditions $v_i+w_i\ge 0$ and the commutation rules of $D_X$. 
In particular, $F^{(\vv,\w)}_0{D_X}$ is a sub-ring of $D_X$, and $F^{(\vv,\w)}_m{D_X}$'s are $F^{(\vv,\w)}_0{D_X}$-submodules of $D_X$. 
We can define the associated graded ring of $D_X$ with respect to the filtration 
\[
\gr^{(\vv,\w)}{D_X}=\bigoplus_{m\in \Z} F^{(\vv,\w)}_m{D_X}/ F^{(\vv,\w)}_{m-1}{D_X}. 
\]
\begin{Definition}
The {\it order} of $P\in D_X$ is defined by 
\[
\ord_{(\vv,\w)}(P)=\min\{m\mid P\in F^{(\vv,\w)}_m{D_X}\}. 
\]
For a non-zero $P\in D_X$ with $\ord_{(\vv,\w)}(P)=m$, the {\it initial form} $\initial_{(\vv,\w)}(P)$ of $P$ is the image of $P$ in $\gr_m^{(\vv,\w)}{D_X}:=F^{(\vv,\w)}_m{D_X}/ F^{(\vv,\w)}_{m-1}{D_X}$. 
For a left ideal $I\subset D_X$, the {\it initial ideal} $\initial_{(\vv,\w)}(I)$ of $I$ is the left ideal of 
$\gr^{(\vv,\w)}{D_X}$ generated by all initial forms of elements in $I$. 
A finite subset $G$ of $D_X$ is called {\it Gr\"obner basis} of $I$ with respect to $(\vv,\w)$ if $I$ is generated by $G$ and 
$\initial_{(\vv,\w)}(I)$ is generated by initial forms of elements in $G$. 
\end{Definition}
It is known that there is an algorithm for computing Gr\"obner bases (\cite{SST} Algorithm 1.2.6). 
We can compute the restriction of ideals to sub-algebras using Gr\"obner bases as in the commutative case. 
\begin{Lemma}\label{elimination of variables}
Let $Z$ be a subsystem of $(\x,\partial_\x)$, and $\C\langle Z \rangle$ a sub-algebra of $D_X$ generated by $Z$ over $\C$.  
Let $(\vv,\w)$ be a weight vector such that $v_i>0$ (resp. $w_j>0$) if $x_i$ (resp. $\partial_{x_j}$) is not a member of $Z$, 
and $v_i=0$ (resp. $w_j=0$) otherwise. 
Let $I$ be a left ideal of $D_X$ and $G$ a Gr\"obner basis of with respect to $(\vv,\w)$, 
then $G\cap \C\langle Z\rangle$ is a system of generators of the left ideal $I\cap \C\langle Z\rangle$. 
\end{Lemma}
We can also compute the intersection of ideals using elimination of variables as in the commutative case. 
\begin{Lemma}\label{intersection}
Let $I$ and $J$ be left ideals of $D_X$. Then 
\[
I\cap J=D_X[u](u I+(1-u)J)\cap D_X. 
\]
\end{Lemma}
\begin{proof}
If $P\in I\cap J$, then $P=uP+(1-u)P\in D_X[u](u I+(1-u)J)\cap D_X$. 
Let $P \in D_X[u](u I+(1-u)J)\cap D_X$. By substituting $1$ and $0$ to $u$, we have $P\in I$ and $P\in J$. 
\end{proof}
Note that substituting some $p\in D_X$ to the variable $u$ makes sense only when $p$ is in the center of $D_X$, that is, $p\in\C$. 
In this case, the left ideal of $D_X[u]$ generated by $u-p$ is a two-side ideal. 

From now on, we assume that the weight vector $(\vv,\w)$ satisfies 
\[
v_i+w_i=0 \mathrm{~for~}i=1,2,\dots,d. 
\]
Then $D_X$ has a structure of a graded algebra:  
We set 
\[
[D_X]^{(\vv,\w)}_m:=\{\sum_{\vv\cdot \boldsymbol{\mu}+\w\cdot \boldsymbol{\nu}= m} 
a_{\boldsymbol{\mu\nu}}\x^{\boldsymbol{\mu}}{\partial_\x}^{\boldsymbol{\nu}}\mid a_{\boldsymbol{\mu\nu}}\in \C\}\subset D_X. 
\]
Then $F^{(\vv,\w)}_m{D_X}=\bigoplus_{k\le m}[D_X]_k^{(\vv,\w)}$ and $\gr_m^{(\vv,\w)}{D_X}\cong [D_X]^{(\vv,\w)}_m$ 
since the commutation rules of $D_X$ are homogeneous of weight $0$. 
Hence $D_X$ is a graded algebra $D_X=\bigoplus_{m\in \Z} [D_X]^{(\vv,\w)}_m$ and isomorphic to $\gr^{(\vv,\w)}{D_X}$. 
In particular $[D_X]^{(\vv,\w)}_0$ is a sub-ring of $D_X$. We call an element in 
$[D_X]^{(\vv,\w)}_m$ a {\it homogeneous element} of degree $m$. 
A left ideal $J$ of $D_X$ is called a {\it homogeneous ideal} if $J$ is generated by homogeneous elements.  

\begin{Definition}
For $P=\sum P_m\in D_X$ with $P_m\in [D_X]^{(\vv,\w)}_m$ and $m_0=\ord_{(\vv,\w)}(P)$, 
we define the homogenization of $P$ with homogenizing variable $u_1$ to be 
\[
P^h=\sum P_mu_1^{m_0-m}\in D_X[u_1]. 
\]
where  For a left ideal $J$ of $D_X$, we define the homogenization of $J$ to be the left ideal of $D_X[u_1]$ 
\[
J^h=\langle P^h\mid P\in J\rangle. 
\]
\end{Definition}
\begin{Definition}
For a left ideal $J$ of $D_X$, we set 
\[
J^*=J^h\cap D_X=\bigoplus_{m\in \Z} (J\cap [D_X]^{(\vv,\w)}_m), 
\]
the ideal of $D_X$ generated by all homogeneous elements in $J$. 
\end{Definition}
\begin{Lemma}\label{homogeneous part}
Let $J=\langle P_1,\dots,P_r\rangle$ be a left ideal of $D_X$. Then 
\[
J^*=D_X[u_1,u_2]\langle P_1^h,\dots,P_r^h,u_1u_2-1\rangle\cap D_X. 
\]
\end{Lemma}
\begin{proof}
It is easy to see that 
\[
J^h=(D_X[u_1,u_1^{-1}]\langle P_1^h,\dots,P_r^h\rangle) \cap D_X[u_1]=\langle P_1^h,\dots,P_r^h,u_1u_2-1\rangle\cap D_X[u_1]. 
\]
Since $J^*=J^h\cap D_X$, we obtain the assertion. 
\end{proof}
\subsection{Bernstein-Sato polynomials}
Budur--Musta\c{t}\u{a}--Saito introduced generalized Bernstein-Sato polynomials (or $b$-function) of arbitrary varieties in \cite{BMS} 
and proved relations between generalized Bernstein-Sato polynomials and multiplier ideals 
using the theory of the $V$-filtration of Kashiwara and Malgralge. 

Let $X$ be the affine space $\C^{n}$ with the coordinate ring $\C[\x]=\C[x_1,\dots,x_n]$, and fix 
an ideal $\a$ of $\C[\x]$ with a system of generators $\f =(f_1,\dots,f_r)$. 
Let $Y=X\times \C^r$ be the affine space $\C^{n+r}$ with the coordinate system $(\x,\t) = (x_1,\dots, x_n,t_1,\dots, t_r)$.  
Then $X\times\{0\}=V(t_1,\dots,t_r)\cong X$ is a linear subspace of $Y$ with the defining ideal 
$I_{X\times\{0\}}=\langle t_1,\dots,t_r\rangle$. 
We denote the rings of differential operators of $X$ and $Y$ by 
\begin{eqnarray*}
D_X&=&\C\langle \x, {\partial_\x}\rangle=\C\langle x_1,\dots,x_n,\partial_{x_1},\dots,\partial_{x_n}\rangle, \\
D_Y&=&\C\langle \x,\t,{\partial_\x},{\partial_\t}\rangle=\C\langle x_1,\dots, x_n, t_1,\dots, t_r,\partial_{x_1},\dots,\partial_{x_n},
\partial_{t_1},\dots,\partial_{t_r}\rangle.  
\end{eqnarray*}
We use the notation 
$\x^{\boldsymbol{\mu}_1}=\prod_{i=1}^n x_i^{{\mu}_{1i}}$, 
$\t^{\boldsymbol{\mu}_2}=\prod_{j=1}^r t_j^{{\mu}_{2i}}$, 
${\partial_\x}^{\boldsymbol{\nu}_1}=\prod_{i=1}^n \partial_{x_i}^{{\nu}_{1i}}$, 
and ${\partial_t}^{\boldsymbol{\nu}_2}=\prod_{j=1}^r\partial_{t_j}^{{\nu}_{2j}}$ 
for $\boldsymbol{\mu}_1 = ({\mu}_{11}, \dots , {\mu}_{1n})$, $\boldsymbol{\nu}_1=({\nu}_{11},\dots,{\nu}_{1n})\in\Z_{\ge 0}^n$ 
and $\boldsymbol{\mu}_2 = ({\mu}_{21},\dots , {\mu}_{2r})$, $\boldsymbol{\nu}_2=({\nu}_{21},\dots,{\nu}_{2r}) \in \Z_{\ge 0}^r$. 
The $\C[\x]$-module $N_{\f}:=\C[\x][\prod_i f_i^{-1},s_1,\dots,s_r]\prod_i f_i^{s_i}$, 
where $s_i$'s are independent variables and $\prod_i f_i^{s_i}$ is a symbol, has a $D_X$-module structure as follows:  
The action of $\C[\x]$ on $N_{\f}$ is given by the canonical one, and the action of $\partial_{x_i}$ are given by
\[
\partial_{x_j} (h\prod_i f_i^{s_i})=\biggl(\partial_{x_j}(h)+h\sum_{k=1}^{r} s_j\frac{\partial_{x_j}(f_k)}{f_k}\biggr)\prod_i f_i^{s_i}.
\] 
for $h\in \C[\x][\prod_i f_i^{-1},s_1,\dots,s_r]$. 
This action is defined formally, but it has an obvious meaning when some integers are substituted for $s_i$'s. 
We define $D_X$-linear actions $t_j$ and $\partial_{t_j}$ on $N_{\f}$ by 
\[
t_j(h(x,s_1,\dots,s_r)\prod_i f_i^{s_i})=h(x,s_1,\dots,s_j+1,\dots,s_r)f_j\prod_i f_i^{s_i}
\] 
and 
\[
\partial_{t_j}(h(x,s_1,\dots,s_r)\prod_i f_i^{s_i})=-s_jh(x,s_1,\dots,s_j-1,\dots,s_r)f_j^{-1}\prod_i f_i^{s_i}. 
\]
for $h(x,s_1,\dots,s_r)\in \C[\x][\prod_i f_i^{-1},s_1,\dots,s_r]$. 
Then it follows that $N_{\f}$ is a $D_Y$-module because the actions defined above respect the commutation rules of $D_Y$. 
Note that $-\partial_{t_i}{t_i} \prod_i f_i^{s_i}=s_i\prod_i f_i^{s_i}$ for all $i$. 
\begin{Definition}[\cite{BMS}]\label{def of b-functions}
Let $\sigma =-(\sum_i \partial_{t_i}{t_i})$, and let $s$ be a new variable. 
Then the global generalized Bernstein-Sato polynomial $b_{\a,g}(s)\in \C [s]$ of $\a=\langle f_1,\dots,f_r\rangle$ and 
$g\in \C[\x]$ is defined to be the monic polynomial of minimal degree satisfying 
\begin{equation}\label{exp1}
b(\sigma)g\prod_i f_i^{s_i}=\sum_{j=1}^{r}P_jgf_j\prod_i f_i^{s_i}
\end{equation}
for some $P_j\in D_X\langle -\partial_{t_j}t_k \mid  1\le j,k \le r\rangle$. 
We define $b_\a(s)=b_{\a,1}(s)$. 

For a prime ideal $\p$ of $\C[\x]$, we define the local generalized Bernstein-Sato polynomial $b^\p_{\a,g}(s)$ at $\p$ 
to be the monic polynomial of minimal degree satisfying 
\begin{equation*}
b(\sigma)gh\prod_i f_i^{s_i}=\sum_{j=1}^{r}P_jgf_j\prod_i f_i^{s_i}
\end{equation*}
for some $P_j\in D_X\langle -\partial_{t_j}t_k \mid  1\le j,k \le r\rangle$ and $h \not\in \p$. 
We define $b^\p_\a(s)=b^\p_{\a,1}(s)$. 
\end{Definition}
Note that $\C[x]_\p \otimes_{\C[x]} D_X$ is the ring of differential operators of $\Spec \C[\x]_\p$. 
It is proved in \cite{BMS} that generalized Bernstein-Sato polynomials are well-defined, 
that is, they do not depend on the choice of generators of $\a$, and all their roots are negative rational numbers. 
These facts follow from the theory of $V$-filtrations of Kashiwara \cite{Ka} and Malgrange \cite{M}. 
When $\a$ is a principal ideal generated by $f$, then $b_\a(s)$ coincides with the classical Bernstein-Sato polynomial $b_f(s)$ of $f$. 
\subsection{V-filtrations}
We will briefly recall the definition and some basic properties of $V$-filtrations. 
See \cite{M}, \cite{Ka}, \cite{Sab} and \cite{BMS} for details. 

We fix the weight vector $(\w,-\w)\in \Z^{2(n+r)}$, $\w=((0,\dots,0),(1,\dots,1))\in \Z^{n}\times\Z^{r}$, that is, 
we assign the weight $1$ to $\partial_{t_j}$, $-1$ to $t_j$, and $0$ to $x_i$ and $\partial_{x_i}$. 
Then 
\[
F_{m}^{(\w,-\w)} D_Y=\{\sum_{-|\boldsymbol{\mu}_2|+|\boldsymbol{\nu}_2|\le m} 
a_{\boldsymbol{\mu}_1 \boldsymbol{\mu}_2 \boldsymbol{\nu}_1 \boldsymbol{\nu}_2}\x^{\boldsymbol{\mu}_1}\t^{\boldsymbol{\mu}_2}{\partial_\x}^{\boldsymbol{\nu}_1}{\partial_\t}^{\boldsymbol{\nu}_2}
\mid a_{\boldsymbol{\mu}_1 \boldsymbol{\mu}_2 \boldsymbol{\nu}_1 \boldsymbol{\nu}_2}\in \C\}
\]
In this paper, we call the decreasing filtration $V^{m}D_Y:=F_{-m}^{(\w,-\w)} D_Y$ on $D_Y$ 
the {\it V-filtration} of $D_Y$ along $X\times\{0\}$ 
(some author call the increasing filtration $F^{(\w,-\w)}$ the $V$-filtration). 
Note that 
\begin{eqnarray*}
V^m{D_Y}&=&\{\sum_{|\boldsymbol{\mu}_2|-|\boldsymbol{\nu}_2|\ge m} 
a_{\boldsymbol{\mu}_1 \boldsymbol{\mu}_2 \boldsymbol{\nu}_1 \boldsymbol{\nu}_2}\x^{\boldsymbol{\mu}_1}\t^{\boldsymbol{\mu}_2}{\partial_\x}^{\boldsymbol{\nu}_1}{\partial_\t}^{\boldsymbol{\nu}_2}\mid 
a_{\boldsymbol{\mu}_1 \boldsymbol{\mu}_2 \boldsymbol{\nu}_1 \boldsymbol{\nu}_2}\in \C\}\\
&=&\{P \in D_Y \mid  P(I_{X\times\{0\}} )^j \subset (I_{X\times\{0\}} )^{j+m} ~\mathrm{for~~any}~ j \ge 0\},
\end{eqnarray*}
with the convention $I_{X\times\{0\}}^{j}=\C[\x,\t]$ for all $j\le 0$. 
\begin{Definition}\label{Vfiltration}
The $V$-filtration along $X\times\{0\}$ on a finitely generated left $D_Y$-module $M$ is an exhaustive decreasing filtration 
$\{V^{{\alpha}}M\}_{\alpha\in \Q}$ indexed by $\Q$, such that: 
\\ {\rm(i)} $V^{{\alpha}}M$ are finitely generated $V^0D_Y$-submodules of $M$. 
\\ {\rm(ii)} 
$\{V^{{\alpha}}M\}_{\alpha}$ is left-continuous and discrete, that is, 
$V^{\alpha} M=\bigcap_{{\alpha}'<{\alpha}}V^{{\alpha}'}M$, and every interval contains only finitely many ${\alpha}\in \Q$ 
with $\gr_V^{\alpha} M\neq 0$. 
Here $\gr_V^{\alpha} M := V^{{\alpha}} M/(\bigcup_{{\alpha}'>{\alpha}} V^{{\alpha}'}M)$. 
\\ {\rm(iii)} 
$(V^iD_Y)(V ^{{\alpha}}M) \subset V^{{\alpha}+i}M$ for any $i \in \Z$, ${{\alpha}} \in \Q$.
\\ {\rm(iv)} 
$(V^iD_Y)(V ^{{\alpha}}M) = V^{{\alpha}+i}M$ for any $i > 0$ if ${{\alpha}}\gg 0$.
\\ {\rm(v)} 
the action of $\sigma+{\alpha}$ is nilpotent on $\gr_V^{\alpha} M$. 
\end{Definition}
\begin{Remark}\label{b from V}
(i) The filtration $V$ is unique if it exists (\cite{Ka}), 
and $D_Y$-submodule $D_Y \prod_i f_i^{s_i}\cong D_Y/\Ann_{D_Y}\prod_i f_i^{s_i}$ of $N_{\f}$ has such a $V$-filtration (see \cite{BMS}). \\
(ii) For ${\alpha}\neq z\in \C$, the action of $\sigma+z$ on  $\gr_V^{\alpha} M$ is invertible. 
Hence, if $\gr_V^{\alpha} M\neq 0$, $u \not\in V^{{\alpha} +\varepsilon}M$ and $b(\sigma)u \in V^{{\alpha} +\varepsilon}M$ 
for some $b(s)\in \C[s]$ and all sufficiently small $\varepsilon >0$, then $s+{\alpha}$ is a factor of $b(s)$. 
\end{Remark}
Let $\iota : X \rightarrow Y$ be the graph embedding $x \mapsto (x, f_1(x),\dots, f_r(x))$ of $\f =(f_1,\dots,f_r)$, 
and $M_{\f}=\iota_+\C[\x]$, where $\iota_+$ denotes the direct image for left $D$-modules.  
There is a natural isomorphism $M_{\f} \cong \C[\x]\otimes_{\C}\C[\partial_{t_1},\dots, \partial_{t_r}]$ (see \cite{Bo}), 
and the action of $\C[\x]$ and $\partial_{t_1}, \dots, \partial_{t_r}$ on $M_{\f}$ is given by the canonical one, 
and the action of a vector field $\xi$ on $X$ and $t_j$ are given by
\begin{eqnarray*}
\xi (g \otimes {\partial_t}^{\nu} ) &=& \xi g \otimes {\partial_t}^{\nu} -\sum_j (\xi f_j)g\otimes \partial_{t_j} {\partial_t}^{{\nu}},\\
t_j(g\otimes {\partial_t}^{\nu}) &=& f_jg \otimes {\partial_t}^{\nu} - {\nu}_j g \otimes {\partial_t}^{{\nu}-1_j}.
\end{eqnarray*}
where $1_j$ is the element of $\Z^r$ whose $i$-th component is $1$ if $i=j$ and $0$ otherwise. 
\begin{Definition}[\cite{BMS}]
Let M be a $D_Y$-module with $V$-filtration. 
For $u\in M$, the Bernstein-Sato polynomial $b_u(s)$ of $u$ is the monic minimal polynomial
of the action of $\sigma$ on $V^0D_Yu/V^1D_Yu$. 
\end{Definition}
By the properties of $V$-filtration in Definition \ref{Vfiltration}, 
the induced filtration $V$ on $(V^0D_Y)u/(V^1D_Y)u$ is finite (see \cite{BMS} Section 2.1) 
This guarantees the existence of $b_{u}(s)$. 
If $u\in V^{\alpha} M$, than $V^0D_Yu \subset V^{\alpha} M$ and $V^1D_Yu \subset V^{{\alpha}+1} M$. 
Hence, if we set ${\alpha}_0=\max\{{\alpha} \mid  u \in V^{\alpha} M\}$, 
then $u\not\in V^{{\alpha}_0+\varepsilon}M$ and $b_u(\sigma)u\in V^1D_Yu \subset V^{{\alpha}_0+1} M\subset V^{{\alpha}_0+\varepsilon}M$ 
for sufficiently small $\varepsilon >0$. Hence 
\[
\max\{{\alpha} \mid  u \in V^{\alpha} M\} = \min\{{\alpha} \mid  \gr_V^{\alpha}((V^0D_Y)u) \neq 0\} 
 \min\{{\alpha} \mid  b_u(-{\alpha}) = 0\}. 
\]
Therefore we conclude the next proposition. 
\begin{Proposition}[\cite{Sab}]
Let M be a $D_Y$-module with $V$-filtration. 
Then 
\[
V^{\alpha} M = \{ u \in M \mid  {\alpha}\le {\alpha}' \mbox{~~if~~} b_u(-{\alpha}') = 0 \}.
\]
\end{Proposition}
Since we have a canonical injection $M_{\f} \rightarrow N_{\f}=\C[\x][\prod_i f_i^{-1},s_1,\dots,s_r]\prod_i f_i^{s_i}$ 
that sends $g\otimes {\partial_t}^{\nu}$ to $g{\partial_t}^{\nu}\prod_i f_i^{s_i}$, 
the generalized Bernstein-Sato polynomial $b_{\a,g}(s)$ coincides with $b_u(s)$ where $u=g\otimes 1 \in M_{\f}$ 
(see Observation \ref{obs} in the next section). 
\subsection{Multiplier ideals}
We will recall the relations between generalized Bernstein-Sato polynomials and multiplier ideals following \cite{BMS}. 
The reader is referred to \cite{L} for general properties of multiplier ideals. 
For a positive rational number $c$, the multiplier ideal $\J(\a^c)$ is defined via a log resolution of $\a$. 
Let $\pi: \tilde{X}\to X=\Spec \C[\x]$ be a log resolution of $\a$, 
namely, $\pi$ is a proper birational morphism, $\tilde{X}$ is smooth, 
and there exists an effective divisor $F$ on $\tilde{X}$ such that $\a\sO_{\tilde{X}}=\sO_{\tilde{X}}(-F)$ 
and the union of the support of $F$ and the exceptional divisor of $\pi$ has simple normal crossings. 
For a given real number $c\ge 0$, the multiplier ideal $\J(\a^c)$ associated to $c$ is defined to be the ideal 
\[
\J(\a^c)=H^0(\tilde{X},\sO_{\tilde{X}}( K_{\tilde{X}/X}-\lfloor cF \rfloor))
\]
where $K_{\tilde{X}/X}=K_{\tilde{X}}-\pi^*K_{X}$ is the relative canonical divisor of $\pi$. 
This definition is independent of the choice of a log resolution $\pi:\tilde{X}\to X$. 
The reader is referred to \cite{L} for general properties of multiplier ideals. 
By the definition, if $c<c'$, then $\J(\a^c)\supset \J(\a^{c'})$ and 
$\J(\a^c)$ is right-continuous in $c$, that is, $\J(\a^c)=\J(\a^{c+\varepsilon})$ for sufficiently small $\varepsilon>0$. 
The multiplier ideals give a decreasing filtration on $\sO_X$, 
and there are rational numbers $0 =c_0< c_1 < c_2 < \cdots$ such that 
$\J (\a^{c_j}) = \J(\a^c) \neq \J(\a^{c_{j+1}})$ for $c_j \leq c< c_{j+1}$. 
These $c_j$ for $j > 0$ are called the {\it jumping coefficients}, 
and the minimal jumping coefficient $c_1$ is called the {\it log-canonical threshold} of $\a$ and denoted by $\lct(\a)$. 
By the definition, it follows that multiplier ideals are integrally closed, 
and the multiplier ideal associated to the log canonical threshold is radical. 
It is known that $\J(\a^{c})=\a\J(\a^{c-1})$ for $c\ge \lambda(\a)$ where $\lambda(\a)$ is the analytic spread of $\a$. 
Recall that analytic spread of $\a$ is the minimal number of elements needed to generate $\a$ up to integral closure, 
and thus $\mu(\a)\ge\lambda(\a)$ where $\mu(\a)$ is the minimal number of generators of $\a$. 
In particular, if $\a$ is a principal ideal generated by $f$, then $\J(f^c)=f\J(f^{c-1})$ for $c\ge 1$. 

Budur--Musta\c{t}\u{a}--Saito proved that the $V$-filtration on $\C[x]$ is essentially equivalent to the filtration by multiplier ideals 
using the theory of mixed Hodge modules (\cite{Sa1}, \cite{Sa2}), 
and gave a description of multiplier ideals in terms of generalized Bernstein-Sato polynomials. 
\begin{Theorem}[\cite{BMS}]\label{multiplier ideals and b-functions}
We denote by $V$ the filtration on $\C[\x]\cong \C[\x]\otimes 1$ induced by the $V$-filtration on $\iota_+ \C[\x]$. 
Then $\J(\a^c) = V^{c+\varepsilon}\C[\x]$ and $V^{\alpha} \C[\x] = \J (\a^{{\alpha}-\varepsilon})$ 
for any ${\alpha}\in \Q$ and $0<\varepsilon \ll 1$. 
Therefore the following hold: \\
{\rm (i)} For a given rational number $c\ge 0$ and a prime ideal $\p\subset \C[x]$, 
\begin{eqnarray*}
\J(\a^c) &=& \{g\in \C[\x] \mid  c < c' \mbox{~~if~~} b_{\a,g}(-c') = 0\}, \\
\J(\a^c)_\p  \cap \C[\x] &=& \{g\in \C[\x] \mid  c < c' \mbox{~~if~~} b^\p_{\a,g}(-c') = 0\}. 
\end{eqnarray*}
In particular, the log canonical threshold $\lct(\a)$ of $\a=\langle f_1,\dots,f_r\rangle$ is the minimal root of $b_\a(-s)$. \\
{\rm (ii)} All jumping coefficients of $\a$ in $[\lct(\a), \lct(\a) + 1)$ are roots of $b_\a(-s)$. 
\end{Theorem}
Therefore an algorithm for computing generalized Bernstein-Sato polynomials induces an algorithm for 
solving membership problem for multiplier ideals, and in particular, an algorithm for computing log canonical thresholds. 
\section{Algorithms for computing generalized Bernstein-Sato polynomials}
In this section, we obtain an algorithm for computing generalized Bernstein-Sato polynomials of arbitrary ideals. 
The algorithms for computing classical Bernstein-Sato polynomials are given by Oaku (see \cite{O1}, \cite{O2}, \cite{O3}). 
We will generalize Oaku's algorithm to arbitrary $n$ and $g$. 

Let $\C[\x]=\C[x_1,\dots,x_n]$ be a polynomial ring over $\C$, 
and $\a$ an ideal with a system of generators $\f =(f_1,\dots,f_r)$, 
and we fix the weight vector $(\w,-\w)\in \Z^{n+r}\times\Z^{n+r}$, $\w=((0,\dots,0),(1,\dots,1))\in \Z^{n}\times\Z^{r}$. 
\begin{Observation}\label{obs}
We rewrite the definition of generalized Bernstein-Sato polynomials in several ways. 

Recall that $b_{\a,g}(s)$ is the monic polynomial of minimal degree satisfying 
\[
b(\sigma)g\prod_i f_i^{s_i}=\sum_{j=1}^{r}P_jgf_j\prod_i f_i^{s_i}
\]
for some $P_j\in D_X\langle -\partial_{t_j}t_k \mid  1\le j,k \le r\rangle$ (Definition \ref{def of b-functions} (\ref{exp1}))
Since 
\begin{eqnarray*}
[D_Y]^{(\w,-\w)}_0&=&\{\sum_{|\boldsymbol{\nu}_2|-|\boldsymbol{\mu}_2|=0} 
a_{\boldsymbol{\mu}_1 \boldsymbol{\mu}_2 \boldsymbol{\nu}_1 \boldsymbol{\nu}_2}\x^{\boldsymbol{\mu}_1}\t^{\boldsymbol{\mu}_2}{\partial_\x}^{\boldsymbol{\nu}_1}{\partial_\t}^{\boldsymbol{\nu}_2}
\mid a_{\boldsymbol{\mu}_1 \boldsymbol{\mu}_2 \boldsymbol{\nu}_1 \boldsymbol{\nu}_2}\in \C\}\\
&=& D_X\langle -\partial_{t_j}t_k \mid  1\le j,k \le r\rangle, 
\end{eqnarray*} 
and $\sigma =-\sum \partial_{t_i}t_i$ is a homogeneous element of degree $0$,
the condition (\ref{exp1}) is equivalent to saying that 
\begin{eqnarray*}
&&(b(\sigma)g-\sum_{j=1}^{r}P_jgf_j)\prod_i f_i^{s_i}=0 \mbox{\quad for~~} {}^\exists P_j\in [D_Y]^{(\w,-\w)}_0 \\
&\Longleftrightarrow&
b(\sigma)g-\sum_{j=1}^{r}P_jgf_j\in \Ann_{[D_Y]^{(\w,-\w)}_0}\prod_i f_i^{s_i} \mbox{\quad for~~} {}^\exists P_j\in [D_Y]^{(\w,-\w)}_0 
\nonumber\\
&\Longleftrightarrow&
b(\sigma)g\in \Ann_{[D_Y]^{(\w,-\w)}_0}\prod_i f_i^{s_i}+[D_Y]^{(\w,-\w)}_0 g\a. 
\end{eqnarray*}
Since $(I+J)\cap [D_Y]^{(\w,-\w)}_0=I\cap [D_Y]^{(\w,-\w)}_0+J\cap [D_Y]^{(\w,-\w)}_0$ for homogeneous ideals $I$ and $J$, 
we have 
\begin{eqnarray*}
&&\Ann_{[D_Y]^{(\w,-\w)}_0}\prod_i f_i^{s_i}+[D_Y]^{(\w,-\w)}_0 g\a\\
=&& (\Ann_{D_Y}\prod_i f_i^{s_i})^*\cap [D_Y]^{(\w,-\w)}_0+(D_Y g\a)\cap[D_Y]^{(\w,-\w)}_0 \\
=&&((\Ann_{[D_Y]}\prod_i f_i^{s_i})^*+D_Y g\a)\cap [D_Y]^{(\w,-\w)}_0. \\
\end{eqnarray*}
Hence the condition (\ref{exp1}) is equivalent to saying that 
\begin{equation}\label{exp2}
b(\sigma)g\in (\Ann_{D_Y}\prod_i f_i^{s_i})^*+D_Y g\a. 
\end{equation}
Since $t_j\prod_i f_i^{s_i}=f_j\prod_i f_i^{s_i}$, the condition (\ref{exp1}) is also equivalent to saying that 
\begin{eqnarray}
\label{exp3}&&b(\sigma)g\prod_i f_i^{s_i}\in  (V^{1}{D_Y})g\prod_i f_i^{s_i}=(F^{(\w,-\w)}_{-1}{D_Y})g\prod_i f_i^{s_i} \\
\label{exp4}\Longleftrightarrow && b(\sigma) \in \initial_{(-w,w)}(\Ann_{D_Y} g\prod_i f_i^{s_i}) \\
\label{exp5}\Longleftrightarrow && b(\sigma)g \in \initial_{(-w,w)}((\Ann_{D_Y} \prod_i f_i^{s_i})\cap D_Yg).   
\end{eqnarray}
By the expression (\ref{exp3}), the generalized Bernstein-Sato polynomial $b_{\a,g}(s)$ coincides with $b_u(s)$ where $u=g\otimes 1 \in M_{\f}$. 
\end{Observation}
By (\ref{exp4}), the polynomial $b_{\a,g}(-s-r)$ coincides with the $b$-function for $\Ann_{D_Y}g\prod_i f_i^{s_i}$ 
with the weight vector $(\w,-\w)$ in \cite{SST}, p.194. 
In the case $g=1$, one can compute $b_\a(s)$ using loc. cit., p.196, Algorithm 5.1.6 by the next lemma. 
\begin{Lemma}\label{annfs}
\[
\Ann_{D_Y} \prod_i f_i^{s_i}=\langle t_i-f_i\mid 1\le i \le r\rangle + \langle \partial_{x_j}
+\sum_{i=1}^{r} \partial_{x_j}(f_i)\partial_{t_i}\mid 1\le j \le n\rangle. 
\]
\end{Lemma}
\begin{proof}
One can prove the assertion similarly to the case $n=1$. See \cite{SST} Lemma 5.3.3. 
\end{proof}
By this lemma and Lemma \ref{homogeneous part}, the homogeneous left ideals 
\[
(\Ann_{D_Y}\prod_i f_i^{s_i})^*+D_Y g\a,\mbox{~~and~~} \initial_{(-w,w)}((\Ann_{D_Y} \prod_i f_i^{s_i})\cap D_Yg),
\] 
in (\ref{exp2}) and (\ref{exp5}) are computable. 
Therefore we can calculate generalized Bernstein-Sato polynomials if we obtain an algorithm for computing the ideal 
$\{b(\x,s)\in \C[\x,s]\mid b(\x,\sigma)\in J \}\cong J \cap \C[x,\sigma]$ for a given homogeneous ideal $J\subset D_Y$. 
One can compute this in the same way as \cite{SST} Algorithm 5.1.6. 
The algorithm calculates $J'=J\cap \C[x,\sigma_1,\dots,\sigma_r]$ first where $\sigma_i=-\partial_{t_i}t_i$, 
then computes $J' \cap \C[x,\sigma]$. This algorithm requires $2r$ new variables. 
We will give an algorithm for computing $J \cap \C[x,\sigma]$ without computing $J'$. 
\begin{Lemma}\label{substitution}
Let $J$ be a homogeneous left ideal of $D_Y$. The following hold: \\
{\rm(i)} $\sigma $ is in the center of $[D_Y]^{(\w,-\w)}_0$.\\
{\rm(ii)} $D_Y[s](J+\langle s-\sigma \rangle)\cap [D_Y]^{(\w,-\w)}_0=J\cap [D_Y]^{(\w,-\w)}_0$.\\
{\rm(iii)} $\{b(\x,s)\in \C[\x,s]\mid b(\x,\sigma)\in J \}=D_Y[s](J+\langle s-\sigma \rangle)\cap \C[\x,s]$. 
\end{Lemma}
\begin{proof}
(i) 
Since the ring $[D_Y]^{(\w,-\w)}_0$ is generated by $\partial_{t_j}{t_k}$, $1\le j,k\le r$, over $D_X$, 
and $\sigma $ commutes with any element of $D_X$, 
it is enough to show that $\sigma (\partial_{t_j}{t_k})=(\partial_{t_j}{t_k})\sigma$ for all $1\le j,k\le r$. 

In the case $j\neq k$, we obtain   
\begin{eqnarray*}
(\partial_{t_j}{t_k})(\partial_{t_j}{t_j})&=&\partial_{t_j}^2t_jt_k,\quad \quad ~~~
(\partial_{t_j}{t_j})(\partial_{t_j}{t_k}) = \partial_{t_j}^2t_jt_k-\partial_{t_j}t_k,\\
(\partial_{t_j}{t_k})(\partial_{t_k}{t_k})&=&\partial_{t_j}\partial_{t_k}t_k^2-\partial_{t_j}t_k,~~
(\partial_{t_k}{t_k})(\partial_{t_j}{t_k}) = \partial_{t_j}\partial_{t_k}t_k^2,
\end{eqnarray*}
thus $(\partial_{t_j}{t_k})(\partial_{t_j}{t_j}+\partial_{t_k}{t_k})=(\partial_{t_j}{t_j}+\partial_{t_k}{t_k})(\partial_{t_j}{t_k})$. 
Hence 
\begin{eqnarray*}
(\partial_{t_j}{t_k})\sigma
&=&-(\partial_{t_j}{t_k})(\partial_{t_j}{t_j}+\partial_{t_k}{t_k}+\sum_{\ell\neq j,k}\partial_{t_j\ell}{t_\ell})\\
&=&-(\partial_{t_j}{t_j}+\partial_{t_k}{t_k})(\partial_{t_j}{t_k})-(\sum_{\ell\neq j,k}\partial_{t_j\ell}{t_\ell})(\partial_{t_j}{t_k})
=\sigma(\partial_{t_j}{t_k}).
\end{eqnarray*}
In the case $j=k$, it is obvious that $(\partial_{t_j}{t_j})\sigma=\sigma(\partial_{t_j}{t_j})$. 
Therefore $\sigma$ is in the center of $[D_Y]^{(\w,-\w)}_0$. \\
(ii) The inclusion $D_Y[s](J+\langle s-\sigma \rangle)\cap [D_Y]^{(\w,-\w)}_0\supset J\cap [D_Y]^{(\w,-\w)}_0$ is trivial. 
We will show the converse inclusion. Let 
\[
h=\sum P_{\ell}s^{\ell}+Q(s)(s-\sigma)\in D_Y[s](J+\langle s-\sigma \rangle)\cap [D_Y]^{(\w,-\w)}_0
\]
where $P_{\ell}\in J$ and $Q(s) \in D_Y[s]$. 
Taking the degree zero part, we may assume $P_{\ell}\in J\cap [D_Y]^{(\w,-\w)}_0$ and $Q[s]\in [D_Y]^{(\w,-\w)}_0[s]$. 
As $\sum P_{\ell}s^{\ell}-\sum P_{\ell}\sigma^{\ell}\in [D_Y]^{(\w,-\w)}_0[s](s-\sigma)$, 
there exists $Q^{'}(s)\in [D_Y]^{(\w,-\w)}_0[s]$ such that 
$h=\sum P_{\ell}\sigma^{\ell}+Q^{'}(s)(s-\sigma)$. 
Since $h\in D_Y$, we have $Q^{'}(s)=0$. Therefore $h=\sum P_{\ell}\sigma^{\ell}=\sum \sigma^{\ell}P_{\ell}\in J\cap [D_Y]^{(\w,-\w)}_0$. \\
(iii) Let $b(\x,\sigma)\in J$. 
Since $b(\x,s)-b(\x,\sigma)\in \langle s-\sigma \rangle$, we have  
\[
b(\x,s)=b(\x,\sigma)+(b(\x,s)-b(\x,\sigma))\in D_Y[s](J+\langle s-\sigma \rangle). 
\] 
Conversely, if $b(\x,s)\in D_Y[s](J+\langle s-\sigma \rangle)\cap \C[\x,s]$, then 
\[
b(\x,\sigma)=b(\x,s)-(b(\x,s)-b(\x,\sigma))\in D_Y[s](J+\langle s-\sigma \rangle). 
\]
Since $b(\x,\sigma)\in [D_Y]^{(\w,-\w)}_0$, and by (ii) , we conclude 
\[
b(\x,\sigma)\in D_Y[s](J+\langle s-\sigma \rangle)\cap [D_Y]^{(\w,-\w)}_0=J\cap [D_Y]^{(\w,-\w)}_0\subset J. 
\]
\end{proof}
\begin{Theorem}[Algorithm for global generalized Bernstein-Sato polynomials 1]\label{Algorithm for global b-functions 1}
Let 
\[
I_{\f}=\langle t_iu_1-f_i\mid 1\le i \le r\rangle + \langle u_1\partial_{x_j}
+\sum_{i=1}^{r} \partial_{x_j}(f_i)\partial_{t_i}\mid 1\le j \le d\rangle+\langle u_1u_2-1\rangle
\]
be a left ideal of $D_Y[u_1,u_2]$. Then compute the following ideals; \\
1. $I_{\f,1}=I_{\f}\cap D_Y$, \\
2. $I_{(\f;g),2}=D_Y[s](I_{\f,1}+g\a+\langle s-\sigma \rangle)\ \cap \C[\x,s]$, \\
3. $I_{(\f;g),3}=I_{(\f;g),2}:g=(I_{(\f;g),2}\cap \langle g \rangle)\cdot g^{-1}$, \\
4. $I_{(\f;g),4}=I_{(\f;g),3}\cap \C[s]$. \\
Then $b_{\a,g}(s)$ is the generator of $I_{(\f;g),4}$. 
\end{Theorem}
\begin{proof}
By Lemma \ref{annfs} and Lemma \ref{homogeneous part}, $I_{\f,1}= (\Ann_{D_Y}\prod_i f_i^{s_i})^*$. 
As $I_{\f,1}+D_Y g\a$ is a homogeneous ideal, we have 
\[
I_{(\f;g),2}=\{b(\x,s)\in \C[\x,s]\mid b(\x,\sigma)\in (\Ann_{D_Y}\prod_i f_i^{s_i})^*+D_Y g\a\}
\]
by Lemma \ref{substitution}. Since $b_{\a,g}(s)$ is the minimal generator of the ideal 
\[
\{b(s)\in \C[s]\mid b(\sigma)g(\x)\in (\Ann_{D_Y}\prod_i f_i^{s_i})^*+D_Y g\a\}, 
\]
it follows that $I_{(\f;g),4}=(I_{(\f;g),2}:g)\cap \C[s]=\langle b_{\a,g}(s)\rangle$. 
\end{proof}
\begin{Theorem}[Algorithm for global generalized Bernstein-Sato polynomials 2]\label{Algorithm for global b-functions 2}
Let 
\[
\tilde{I}_{\f}=\langle t_i-f_i\mid 1\le i \le r\rangle + \langle \partial_{x_j}
+\sum_{i=1}^{r} \partial_{x_j}(f_i)\partial_{t_i}\mid 1\le j \le d\rangle\subset D_Y, 
\]
and compute the following ideals; \\
0. $\tilde{I}_{(\f;g),0}=\tilde{I}_{\f}\cap D_Yg=D_Y[u](u\tilde{I}_{\f}+(1-u)g)\cap D_Y$, \\
1. $\tilde{I}_{(\f;g),1}=\initial_{(\w,-\w)}(\tilde{I}_{(\f;g),0})$, \\
2. $\tilde{I}_{(\f;g),2}=(\tilde{I}_{(\f;g),1}+\langle s-\sigma \rangle)\ \cap \C[\x,s]$, \\
3. $\tilde{I}_{(\f;g),3}=\tilde{I}_{(\f;g),2}:g=\tilde{I}_{(\f;g),2}\cdot g^{-1}$, \\
4. $\tilde{I}_{(\f;g),4}=\tilde{I}_{(\f;g),3}\cap \C[s]$. \\
Then $b_{\a,g}(s)$ is the generator of $\tilde{I}_{(\f;g),4}$. 
\end{Theorem}
\begin{proof}
By Lemma \ref{annfs} and Lemma \ref{intersection}, $\tilde{I}_{(\f;g),0}= (\Ann_{D_Y}\prod_i f_i^{s_i})\cap D_Yg$. 
As $\tilde{I}_{(\f;g),1}$ is a homogeneous ideal, we have 
\[
\tilde{I}_{(\f;g),2}=\{b(\x,s)\in \C[\x,s]\mid b(\x,\sigma)\in \initial_{(\w,-\w)}((\Ann_{D_Y}\prod_i f_i^{s_i})\cap D_Y g) \}
\]
by Lemma \ref{substitution}. 
Since $\tilde{I}_{(\f;g),0}\subset D_Yg$, we have $\tilde{I}_{(\f;g),1}\subset \initial_{(\w,-\w)}(D_Yg)=D_Yg$. 
Hence $\tilde{I}_{(\f;g),2}\subset (D_Yg +\langle s-\sigma \rangle)\ \cap \C[\x,s]= \C[\x,s]g$, 
and $\tilde{I}_{(\f;g),2}:g=\tilde{I}_{(\f;g),2}\cdot g^{-1}$. 
Since $b_{\a,g}(s)$ is the minimal generator of the ideal 
\[
\{b(s)\in \C[s]\mid b(\sigma)g(\x)\in (\Ann_{D_Y}\prod_i f_i^{s_i})^*+D_Y g\a\}, 
\]
it follows that $\tilde{I}_{(\f;g),4}=(\tilde{I}_{(\f;g),2}:g)\cap \C[s]=\langle b_{\a,g}(s)\rangle$. 
\end{proof}
\begin{Remark}\label{same algorithm}
Note that 
\[
I_{(\f;g),3}=\tilde{I}_{(\f;g),3}=\{b(\x,s)\in \C[s]\mid b(\x,\sigma)g(\x)\prod_i f_i^{s_i}\in (V^1D_Y) g(x)\prod_i f_i^{s_i}\}, 
\]
and in the case $g=1$, it follows that 
\[
I_{(\f;1),2}=\tilde{I}_{(\f;1),2}=\{b(\x,s)\in \C[\x,s]\mid b(\x,\sigma)\prod_i f_i^{s_i} \in (V^1D_Y)  \prod_i f_i^{s_i}\}. 
\]
\end{Remark}
We can compute local generalized Bernstein-Sato polynomials similarly to the classical case using primary decompositions 
(see \cite{O1}, \cite{O2}, and \cite{O3}). 
\begin{Theorem}[Algorithm for local generalized Bernstein-Sato polynomials]\label{Algorithm for local b-functions}
Let $I_{(\f;g),3}=\tilde{I}_{(\f;g),3}$ be the ideals in Theorem \ref{Algorithm for global b-functions 1} 
and Theorem \ref{Algorithm for global b-functions 2} 
with primary decompositions $I_{(\f;g),3}=\cap_{i=1}^{\ell} \q_i$. 
Then $b^\p_{\a,g}(s)$ is the generator of the ideal
\[
\bigcap_{\q_i \cap \C[\x] \subset \p}\q_i\cap \C[s].
\]
\end{Theorem}
\begin{proof}
We set $R=\C[\x]$. 
By the definition, $b^\p_{\a,g}(s)$ is the generator of the ideal 
\[
\{b(s)\in \C[s]\mid b(\sigma)g(\x)h(\x)\in (\Ann_{D_Y}\prod_i f_i^{s_i})^*+D_Y g\a \mbox{\quad for~~} {}^\exists h\not\in \p \}.  
\]
This ideal equals to  
\begin{eqnarray*}
&&\{b(s)\in \C[s]\mid b(s)g(\x)h(\x)\in I_{(\f;g),2} \mbox{\quad for~~} {}^\exists h\not\in \p \}\\
&=&\{b(s)\in \C[s]\mid b(s)h(\x)\in I_{(\f;g),3} \mbox{\quad for~~} {}^\exists h\not\in \p \}\\
&=&(I_{(\f;g),3}\otimes_{R} R_\p )\cap \C[s]. 
\end{eqnarray*} 
Since $I_{(\f;g),3}\otimes_{R} R_\p =\bigcap_{\q_i \cap R\subset \p}\q_i\otimes_{R} R_\p$, 
and $\q_i\otimes_{R} R_\p=R_\p[s]$ if and only if $\q_i \cap R \subset \p$, 
it follows that $b^\p_{\a,g}(s)$ is the generator of the ideal  
$(\bigcap_{\q_i \cap R \subset \p}\q_i)\cap \C[s]$. 
\end{proof}
The algorithm for computing generalized Bernstein-Sato polynomials enables us to solve the membership problem for multiplier ideals, 
but does not give a system of generators of multiplier ideals. 
We have to compute $b_{\a,g}(s)$ for infinitely many $g$ to obtain a system of generators. 
\section{Algorithms for computing multiplier ideals}
The purpose of this section is to obtain algorithms for computing the system of generators of multiplier ideals. 
To do this, we modify the definition of Budur--Musta\c{t}\v{a}--Saito's Bernstein-Sato polynomial. 

As in the previous section, let $\C[\x]=\C[x_1,\dots,x_n]$ be a polynomial ring over $\C$, 
and $\a$ an ideal with a system of generators $\f =(f_1,\dots,f_r)$, 
and fix the weight vector $(\w,-\w)\in \Z^{n+r}\times\Z^{n+r}$, $\w=((0,\dots,0),(1,\dots,1))\in \Z^{n}\times\Z^{r}$. 
We set $\delta=1\otimes 1\in M_{\f}=\iota_+ \C[\x]$ and $\overline{M}_{\f}^{(m)}=(V^0D_Y) \delta/(V^m D_Y) \delta$. 
The induced filtration $V$ on the $\overline{M}_{\f}^{(m)}$ is finite by the definition of the $V$-filtration 
(Definition \ref{Vfiltration}) as in the case $m=1$. 
For $g\in \C[\x]$, we denote by $\overline{g\otimes 1}$ the image of $g\otimes 1=g\delta$ in $\overline{M}_{\f}^{(m)}$. 
\begin{Definition}\label{new Bernstein-Sato polynomial}
We define $b_{\a,g}^{(m)}(s)$ to be the monic minimal polynomial of the action of $\sigma$ on 
$(V^0D_Y)\overline{g\otimes 1}\subset \overline{M}_{\f}^{(m)}$. 
We define $b^{(m)}_{\a}=b^{(m)}_{{\a,1}}$. 
\end{Definition}
The existence of $b_{\a,g}^{(m)}(s)$ follows from the finiteness of the filtration $V$ on $\overline{M}_{\f}^{(m)}$, 
and the rationality of its roots follows from the rationality of the $V$-filtration. 
Note that $b_{\a,g}^{(m)}(s)=1$ if and only if $g\otimes 1 \in (V^m D_Y) \delta$. 
\begin{Observation}\label{on new b-function}
Since the ring $V^0D_Y$ is generated by $\t=t_1,\dots,t_r$ over $[D_Y]^{(\w,-\w)}_0$, 
$\sigma$ is in the center of $[D_Y]^{(\w,-\w)}_0$, and $t_i\cdot\overline{g\otimes 1}=f_i\cdot\overline{g\otimes 1}$, 
it follows that $b_{\a,g}^{(m)}(s)$ is the monic polynomial $b(s)$ of minimal degree satisfying $b(\sigma)\overline{g\otimes 1}=0$. 
This is equivalent to saying 
\[
b(s)g \prod_i f_i^{s_i}\in V^mD_Y\prod_i f_i^{s_i}. 
\]
Since $V^mD_Y$ is generated by all monomials in $\t=t_1,\dots,t_r$ of degree $m$ as $V^0D_Y$-module, 
and $t_j\prod_i f_i^{s_i}=f_j\prod_i f_i^{s_i}$, 
our Bernstein-Sato polynomial $b_{\a,g}^{(m)}(s)$ is the monic polynomial $b(s)$ of minimal degree satisfying 
\begin{equation}\label{exp6}
b(s)g \prod_i f_i^{s_i}= \sum_j P_j h_j \prod_i f_i^{s_i} \mbox{\quad for~~} {}^\exists P_j\in [D_Y]^{(\w,-\w)}_0, {}^\exists h_j\in \a^m.
\end{equation}
Hence if $g,h\in\C[\x]$ are polynomials such that $g$ divides $h$, then $b_{\a,h}^{(m)}(s)$ is a factor of $b_{\a,g}^{(m)}(s)$. 
In particular, $b_{\a,g}^{(m)}(s)$ is a factor of $b_{\a}^{(m)}(s)$ for all $g\in\C[\x]$. 
\end{Observation}
We obtain a description of multiplier ideals in terms of our Bernstein-Sato polynomials similarly to 
Theorem \ref{multiplier ideals and b-functions}. 
\begin{Theorem}\label{multiplier ideals and b-functions 2}
{\rm (i)} For a given rational number $c<m+\lct(\a)$, 
\[
\J(\a^c)=\{g \in \C[\x] \mid  c < c' \mbox{~~if~~} b_{\a,g}^{(m)}(-c') = 0 \}. 
\]
In particular, the log canonical threshold $\lct(\a)$ of $\a=\langle f_1,\dots,f_r\rangle$ is the minimal root of $b^{(m)}_{1}(-s)$. \\
{\rm (ii)} All jumping coefficients of $\a$ in $[\lct(\a), \lct(\a) + m)$ are roots of $b^{(m)}_{\a}(-s)$. 
\end{Theorem}
\begin{proof}
First note that we have $\J(\a^c) = V^{c+\varepsilon}\C[\x]$ for $0<\varepsilon\ll 1$ by Theorem \ref{multiplier ideals and b-functions}. 
Since $\lct(\a)=\max\{{c} \mid  \delta \in V^{c} M_{\f}\}$, 
we have $(V^m D_Y) \delta \subset V^{m+\lct(\a)}M_{\f}$. 

If $\in\C[\x]$ is a polynomial with $b_{\a,g}^{(m)}(s)=1$, then $g\otimes 1 \in (V^m D_Y) \delta\subset V^{m+\lct(\a)}M_{\f}$. 
Thus $g\in \J(\a^c)$ for all $c<m+\lct(\a)$. 

If $g\not\in \J(\a^{m+\lct(\a)})$, then $g\otimes 1\not\in V^{m+\lct(\a)} M_{\f}$. 
Hence we obtain 
\[
\max\{{c} \mid  g\otimes 1 \in V^{c} M\}  = \min\{{c} \mid  b_{\a,g}^{(m)}(-{c}) = 0\} 
\]
by Remark \ref{b from V}. 
If $g$ is a polynomial with $b^{(m)}_{\a,g}(-s)=1$, then $g\otimes 1 \in (V^m D_Y) \delta\subset V^{m+\lct(\a)}M_{\f}$. 
Hence $g\in V^{m+\lct(\a)} \C[\x]$. 

Therefore, for ${c}<m+\lct(\a)$, we have 
\[
V^{c} \C[\x] = \{ g \in \C[\x] \mid  {c}\le {c}' \mbox{~~if~~} b_{\a,g}^{(m)}(-{c}') = 0 \}. 
\]
Thus 
\[
V^{c+\varepsilon} \C[\x] = \{ g \in \C[\x] \mid  {c}< {c}' \mbox{~~if~~} b_{\a,g}^{(m)}(-{c}') = 0 \} 
\]
for $0<\varepsilon\ll 1$. 
This proves the assertion. 
\end{proof}
Our Bernstein-Sato polynomials do not tell us any information about multiplier ideals $\J(\a^c)$ for $c\ge m+\lct(\a)$. 
Our definition, however, has the advantage in that we need just one module 
$\overline{M}_{\f}^{(m)}$ to compute Bernstein-Sato polynomials $b^{(m)}_{\a,g}(s)$ for all $g$. 
This fact enables us to compute multiplier ideals. 
\begin{Theorem}[Algorithm for multiplier ideals 1]\label{algorithm for computing multiplier ideals 1}
Let 
\[
I_{\f}=\langle t_iu_1-f_i\mid 1\le i \le r\rangle + \langle u_1\partial_{x_j}
+\sum_{i=1}^{r} \partial_{x_j}(f_i)\partial_{t_i}\mid 1\le j \le d\rangle+\langle u_1u_2-1\rangle
\]
be a left ideal of $D_Y[u_1,u_2]$. Then compute the following ideals; \\
1. $I_{\f,1}=I_{\f}\cap D_Y$, \\
2. $J_{\f}(m)=D_Y[s](I_{\f,1}+\a^m+\langle s-\sigma \rangle)\ \cap \C[\x,s]$. \\
Then the following hold: \\
{\rm (i)}  $b_{\a,g}^{(m)}(s)$ is the generator of $(J_{\f}(m):g)\cap \C[s]$. \\
{\rm (ii)} Let $J_{\f}(m)=\cap_{i=1}^{\ell} \q_i$ be a primary decomposition of $J_{\f}(m)$. 
Then, for $1\le i \le \ell$, 
there exists $c(i)$, a root of $b^{(m)}_{\a}(-s)$, such that the generator of $\q_i\cap \C[s]$ is some power of $s+c(i)$, 
and 
\[
\{c(i)\mid 1\le i \le \ell \}=\{c'\mid b^{(m)}_{\a}(-c')=0\}.  
\]
{\rm (iii)} 
For $c<\lct(\a)+m$, 
\[
\J(\a^c)=\bigcap_{i\in\{j|c(j)\le c\}}\q_i \cap \C[\x]. 
\]
\end{Theorem}
\begin{proof}
(i) As we already see in Observation \ref{on new b-function} \ref{exp6}, 
$b_{\a,g}^{(m)}(s)$ is the monic polynomial $b(s)$ of minimal degree satisfying 
\[
b(s)g \prod_i f_i^{s_i}= \sum_j P_j h_j \prod_i f_i^{s_i}
\]
for some $P_j\in [D_Y]^{(\w,-\w)}_0$ and $h_j\in \a^m$. 
This is equivalent to 
\begin{eqnarray*}
&&(b(\sigma)g-\sum_{j}P_jh_j)\prod_i f_i^{s_i}=0 \mbox{\quad for~~}{}^\exists P_j\in [D_Y]^{(\w,-\w)}_0, {}^\exists h_j\in \a^m \\
&\Longleftrightarrow&
b(\sigma)g-\sum_j P_jh_j\in \Ann_{[D_Y]^{(\w,-\w)}_0} \prod_i f_i^{s_i} 
\mbox{\quad for~~} {}^\exists P_j\in [D_Y]^{(\w,-\w)}_0, {}^\exists h_j\in \a^m \\
&\Longleftrightarrow&
b(\sigma)g\in \Ann_{[D_Y]^{(\w,-\w)}_0}\prod_i f_i^{s_i}+[D_Y]^{(\w,-\w)}_0 \a^m \nonumber\\
&&
\qquad\quad\!\!\!= (\Ann_{D_Y}\prod_i f_i^{s_i})^*\cap [D_Y]^{(\w,-\w)}_0+(D_Y \a^m)\cap[D_Y]^{(\w,-\w)}_0 \nonumber\\
&&
\qquad\quad\!\!\!= ((\Ann_{[D_Y]}\prod_i f_i^{s_i})^*+D_Y \a^m)\cap [D_Y]^{(\w,-\w)}_0 \nonumber\\
&\Longleftrightarrow&
b(\sigma)g\in (\Ann_{D_Y}\prod_i f_i^{s_i})^*+D_Y \a^m. 
\end{eqnarray*}
Hence $b^{(m)}_{\a,g}(s)$ is the generator of the ideal 
\[
\{b(s)\in \C[s]\mid b(\sigma)g(\x)\in (\Ann_{D_Y}\prod_i f_i^{s_i})^*+D_Y\a^m\}. 
\]
On the other hand, we have $I_{\f,1}= (\Ann_{D_Y}\prod_i f_i^{s_i})^*$ by Lemma \ref{annfs} and Lemma \ref{homogeneous part}. 
Since $(\Ann_{D_Y}\prod_i f_i^{s_i})^*+D_Y \a^m$ is a homogeneous ideal, we obtain 
\[
J_{\f}(m)=\{b(\x,s)\in \C[\x,s]\mid b(\x,\sigma)\in (\Ann_{D_Y}\prod_i f_i^{s_i})^*+D_Y\a^m \}, 
\]
by Lemma \ref{substitution}. 
Therefore $b^{(m)}_{\a,g}(s)$ is the generator of $(J_\f(m):g)\cap\C[s]$. \\
(ii) Since $b^{(m)}_{\a}(s)$ is the generator of $J_{\f}(m)\cap \C[s]$, 
and $\bigcap \q_i\cap \C[s]$ is a primary decomposition of $J_{\f}(m)\cap \C[s]$, we conclude (ii). 
Note that $\bigcap \q_i\cap \C[s]$ is not an irredundant primary decomposition of $J_{\f}(m)\cap \C[s]$ in general 
even if $\cap_{i=1}^{\ell} \q_i$ is irredundant, 
and it may happen that $c(i)=c(j)$ for $i\neq j$. \\
(iii) Let $J_c$ be the ideal on the right-hand side of the equality. 
If $g\in \J(\a^c)$, then $b^{(m)}_{\a,g}(s)g\in J_{\f}(m)$ and all the roots of $b^{(m)}_{\a,g}(-s)$ are strictly larger than $c$. 
Suppose that $\q_i$ satisfies $c(i) \le c$. 
Since $b^{(m)}_{\a,g}(s)g\in J_{\f}(m)\subset \q_i$, and any power of $b^{(m)}_{\a,g}(s)$ is not in $\q_i$, we have $g\in \q_i$. 
Hence we conclude $g\in J_c$. 

For the converse inclusion, let $g\in J_c$. 
If $c(i) \le c$, then $g\in \q_i$, and thus $\q_i:g =\C[\x,s]$. 
Since $J_{\f}(m):g=\cap_{i=1}^{\ell} (\q_i:g)$ and by (i), all the roots of $b^{(m)}_{\a,g}(-s)$ are also strictly larger than $c$. 
Therefore we have $g\in \J(\a^c)$ by Theorem \ref{multiplier ideals and b-functions 2} (i).  
\end{proof}
\begin{Theorem}[Algorithm for multiplier ideals 2]\label{algorithm for computing multiplier ideals 2}
Let the notation be as in Theorem \ref{algorithm for computing multiplier ideals 1}. 
Compute $b_{\a}^{(m)}(s)$, the generator of $J_\f(m)\cap \C[s]$, 
and let $c_1,\dots,c_\ell$, $c_i<c_{i+1}$, be all the roots of $b_{\a}^{(m)}(-s)$. 
Then, for $c<\lct(\a)+m$, 
\[
\J(\a^c)=(J_\f(m):(\prod_{c< c_i}(s+c_i))^\infty ) \cap \C[\x]. 
\]
\end{Theorem}
\begin{proof}
Let $g\in \J(\a^c)$. Then $b^{(m)}_{\a,g}(s)g\in J_{\f}(m)$ and all the roots of $b^{(m)}_{\a,g}(-s)$ are strictly larger than $c$. 
Since $b^{(m)}_{\a,g}(s)$ is a factor of $b^{(m)}_{\a}(s)$, it follows that $g$ is in the ideal on the right hand side. 
Hence $\J(\a^c)\subset (J_\f(m):(\prod_{c< c_i}(s+c_i))^\infty ) \cap \C[\x]. $

For the converse inclusion, let $g\in (J_\f(m):(\prod_{c< c_i}(s+c_i))^\infty ) \cap \C[\x]$. 
Then there exists a polynomial $b(s)$ such that $b(s)g\in J_\f(m)$ and the set of roots of $b(-s)$ is in $\{c_i\mid c<c_i\}$. 
Since $b(s)\in (J_\f(m):g)\cap\C[s]=\langle b_{\a,g}^{(m)}(s)\rangle$, $b_{\a,g}^{(m)}(s)$ is a factor of $b(s)$. 
Hence all the roots of $b^{(m)}_{\a,g}(-s)$ are strictly larger than $c$. 
Therefore we have $g\in \J(\a^c)$ by Theorem \ref{multiplier ideals and b-functions 2} (i).  
\end{proof}
\begin{Remark}
(i) In the case $m=1$, $J_{\f}(1)$ coincides with $I_{(\f;1),2}$ in Theorem \ref{Algorithm for global b-functions 1} 
and $b_\a(s)=b_{\a}^{(1)}(s)$. 
Since $I_{(\f;1),2}=\tilde{I}_{(\f;1),2}$ by Remark \ref{same algorithm}, one can obtain another algorithm in this case. \\
(ii) Since  $\J(\a^{c})=\a\J(\a^{c-1})$ for $c\ge \lambda(\a)$, 
it is enough to compute $J_{\f}(m)$ for $m$ satisfying $m\ge \lambda(\a)-\lct(\a)$ to obtain all multiplier ideals. 
In particular, if $\a$ is a principal ideal, it is enough to compute $J_\f(1)$. 
\end{Remark}
\section{Examples}
The computations were made using Kan/sm1 \cite{T} and Risa/Asir \cite{N}.  
\begin{Example}
(i) (\cite{SS}) Let $M=(x_{ij})_{ij}$ be the $n\times n$ general matrix, and $f=\det M\in \C[x_{ij}\mid 1\le i,j\le r]$. 
Then 
$b_f(s)=\prod_{i=1}^r(s+i)$, 
and $b_f(\sigma)f^s=f({\partial_\x})f^{s+1}$. 
Here, we use the notation $h({\partial_\x})$ to mean $\sum a_{\mu} {\partial_\x}^{\mu}$ for $h(\x)=\sum a_{\mu} \x^{\mu}$. \\
(ii) Let $\a=I_2\left(
	\begin{array}{ccc}
	x_1 & x_2 & x_3 \\
	x_4 & x_5 & x_6 
	\end{array}
\right)\subset \C[x_1,\dots,x_6]$ 
with a system of generators $(f_1,f_2,f_3)=(x_1x_5-x_2x_4,x_2x_6-x_3x_5,x_3x_4-x_1x_6)$. 
Then 
\[
b_\a(s)=(s+2)(s+3), 
\]
and 
\[
(\sigma+2)(\sigma+3)f_1^{s_1}f_2^{s_2}f_3^{s_3}=f_1({\partial_\x})f_1^{s_1+1}f_2^{s_2}f_3^{s_3}
+f_2({\partial_\x})f_1^{s_1}f_2^{s_2+1}f_3^{s_3}+f_3({\partial_\x})f_1^{s_1}f_2^{s_2}f_3^{s_3+1}. 
\]
\end{Example}
\begin{Example}
(i) Let $f={x}^2+{y}^3\in \C[{x},{y}]$. Then 
\begin{eqnarray*}
b_f(s)      &=&\bigl(s+\frac{5}{6}\bigr)\bigl(s+1\bigr)\bigl(s+\frac{7}{6}\bigr), \\
b_{f,{x}}(s)&=&\bigl(s+1\bigr)\bigl(s+\frac{11}{6}\bigr)\bigl(s+\frac{13}{6}\bigr), \\
b_{f,{y}}(s)&=&\bigl(s+1\bigr)\bigl(s+\frac{7}{6}\bigr)\bigl(s+\frac{11}{6}\bigr). 
\end{eqnarray*}
(ii) Let $\a=\langle {x}^2,{y}^3\rangle\subset \C[{x},{y}]$. Then 
\begin{eqnarray*}
b_\a(s)      &=&\bigl(s+\frac{5}{6}\bigr)\bigl(s+\frac{7}{6}\bigr)\bigl(s+\frac{4}{3}\bigr)\bigl(s+\frac{3}{2}\bigr)\bigl(s+\frac{5}{3}\bigr)\bigl(s+2\bigr), \\
b_{\a,{x}}(s)&=&\bigl(s+\frac{4}{3}\bigr)\bigl(s+\frac{5}{3}\bigr)\bigl(s+\frac{11}{6}\bigr)\bigl(s+2\bigr)\bigl(s+\frac{13}{6}\bigr)\bigl(s+\frac{5}{2}\bigr), \\
b_{\a,{y}}(s)&=&\bigl(s+\frac{7}{6}\bigr)\bigl(s+\frac{3}{2}\bigr)\bigl(s+\frac{5}{3}\bigr)\bigl(s+\frac{11}{6}\bigr)\bigl(s+2\bigr)\bigl(s+\frac{7}{3}\bigr). 
\end{eqnarray*}
(iii) Let $f=x_3{x}^2+x_4{y}^3\in \C[{x},{y},x_3,x_4]$. Then 
\begin{eqnarray*}
b_f(s)      &=&\bigl(s+\frac{5}{6}\bigr)\bigl(s+1\bigr)\bigl(s+\frac{7}{6}\bigr)\bigl(s+\frac{4}{3}\bigr)\bigl(s+\frac{3}{2}\bigr)\bigl(s+\frac{5}{3}\bigr)\bigl(s+2\bigr), \\
b_{f,{x}}(s)&=&\bigl(s+1\bigr)\bigl(s+\frac{4}{3}\bigr)\bigl(s+\frac{5}{3}\bigr)\bigl(s+\frac{11}{6}\bigr)\bigl(s+2\bigr)\bigl(s+\frac{13}{6}\bigr)\bigl(s+\frac{5}{2}\bigr), \\
b_{f,{y}}(s)&=&\bigl(s+1\bigr)\bigl(s+\frac{7}{6}\bigr)\bigl(s+\frac{3}{2}\bigr)\bigl(s+\frac{5}{3}\bigr)\bigl(s+\frac{11}{6}\bigr)\bigl(s+2\bigr)\bigl(s+\frac{7}{3}\bigr). 
\end{eqnarray*}
\end{Example}
\begin{Example}
Let $f=(x+y)^2-(x-y)^5\in \C[x,y]$ which is not a non-degenerate polynomial. Then 
\begin{eqnarray*}
b_f(s)&=&\bigl(s+\frac{7}{10}\bigr)\bigl(s+\frac{9}{10}\bigr)\bigl(s+ 1 \bigr)\bigl(s+\frac{11}{10}\bigr)\bigl(s+\frac{13}{10}\bigr), \\
b_{f,x}(s)=b_{f,y}(s)&=&\bigl(s+\frac{9}{10}\bigr)\bigl(s+ 1 \bigr)\bigl(s+\frac{11}{10}\bigr)\bigl(s+\frac{13}{10}\bigr)\bigl(s+\frac{17}{10}\bigr), \\
b_{f,x+y}(s)&=&\bigl(s+ 1 \bigr)\bigl(s+\frac{17}{10}\bigr)\bigl(s+\frac{19}{10}\bigr)\bigl(s+\frac{21}{10}\bigr)\bigl(s+\frac{23}{10}\bigr), \\
b_{f,xy}(s) &=&\bigl(s+ 1 \bigr)\bigl(s+\frac{11}{10}\bigr)\bigl(s+\frac{13}{10}\bigr)\bigl(s+\frac{17}{10}\bigr)\bigl(s+\frac{19}{10}\bigr), 
\end{eqnarray*}
and the multiplier ideals are 
\begin{equation*}
\J(f^c)=
\left\{
	\begin{array}{lll}
	\C[x,y] & \mbox{ $0\le c<\frac{7}{10}$,} \\[1mm]
	\langle x, y\rangle=\langle x+y, x-y\rangle & \mbox{ $\frac{7}{10}\le c<\frac{9}{10}$,}\\[1mm]
	\langle x+y, xy\rangle=\langle x+y, (x-y)^2\rangle & \mbox{ $\frac{9}{10}\le c<1$,}
	\end{array}
\right.
\end{equation*}
and $\J(f^c)=f\J(f^{c-1})$ for $c\ge 1$. 
\end{Example}
\begin{Example}
Let $f=xy(x+y)(x+2y)\rangle \subset \C[x,y]$. Then 
\begin{eqnarray*}
b_f(s) &=&\bigl(s+\frac{1}{2}\bigr)\bigl(s+\frac{3}{4}\bigr)\bigl(s+ 1 \bigr)^2\bigl(s+\frac{5}{4}\bigr)\bigl(s+\frac{3}{2}\bigr), \\
b^{(2)}_{f}(s)
&=&\bigl(s+\frac{1}{2}\bigr)\bigl(s+\frac{3}{4}\bigr)\bigl(s+ 1 \bigr)^2\bigl(s+\frac{5}{4}\bigr)\bigl(s+\frac{3}{2}\bigr)
\bigl(s+\frac{7}{4}\bigr)\bigl(s+ 2 \bigr)^2 \\
&&\bigl(s+\frac{9}{4}\bigr)\bigl(s+\frac{5}{2}\bigr), \\
b_{f,x}(s)=b_{f,y}(s)&=&\bigl(s+\frac{3}{4}\bigr)\bigl(s+ 1 \bigr)^2\bigl(s+\frac{5}{4}\bigr)\bigl(s+\frac{3}{2}\bigr)\bigl(s+\frac{7}{4}\bigr)\bigl(s+ 2 \bigr), \\
b_{f,x^2}(s)=b_{f,y^2}(s)&=&\bigl(s+ 1 \bigr)^2\bigl(s+\frac{5}{4}\bigr)\bigl(s+\frac{3}{2}\bigr)\bigl(s+\frac{7}{4}\bigr)\bigl(s+2)\bigl(s+3\bigr), 
\end{eqnarray*}
and the multiplier ideals are 
\begin{equation*}
\J(f^c)=
\left\{
	\begin{array}{lll}
	\C[x,y] & \mbox{ $0\le c<\frac{1}{2}$,} \\[1mm]
	\langle x, y\rangle & \mbox{ $\frac{1}{2}\le c<\frac{3}{4}$,} \\[1mm]
	\langle x, y\rangle^2 & \mbox{ $\frac{3}{4}\le c<1$,}
	\end{array}
\right.
\end{equation*}
and $\J(f^c)=f\J(f^{c-1})$ for $c\ge 1$. 
Note that $\frac{5}{4}$ and $\frac{9}{4}$ are roots of $b^{(2)}_{f}(-s)$ which are not jumping coefficients.  
\end{Example}
\begin{Example}
Let $\a\subset \C[x_1,x_2,x_3]$ be the defining ideal of the space monomial curve $\Spec\C[T^4,T^5,T^6]\subset \C^3$ with a system of generators $\f=(x_2^2-x_1x_3, x_1^3-x_3^2 )$. 
Then generalized Bernstein-Sato polynomials are 
\begin{eqnarray*}
b_\a(s)&=&\bigl(s+\frac{17}{12}\bigr)\bigl(s+\frac{3}{2}\bigr)\bigl(s+\frac{19}{12}\bigr)\bigl(s+\frac{7}{4}\bigr)\bigl(s+\frac{11}{6}\bigr)\bigl(s+\frac{23}{12}\bigr)\bigl(s+2\bigr)\bigl(s+\frac{25}{12}\bigr)\\
&&\bigl(s+\frac{13}{6}\bigr)\bigl(s+\frac{9}{4}\bigr), \\
b_{\a,x_1}(s)&=&\bigl(s+\frac{7}{4}\bigr)\bigl(s+\frac{23}{12}\bigr)\bigl(s+2\bigr)\bigl(s+\frac{25}{12}\bigr)\bigl(s+\frac{13}{6}\bigr)\bigl(s+\frac{9}{4}\bigr)\bigl(s+\frac{29}{12}\bigr)\bigl(s+\frac{5}{2}\bigr)\\
&&\bigl(s+\frac{31}{12}\bigr)\bigl(s+\frac{17}{6}\bigr), \\
b_{\a,x_2}(s)&=&\bigl(s+\frac{11}{6}\bigr)\bigl(s+2\bigr)\bigl(s+\frac{13}{6}\bigr)\bigl(s+\frac{9}{4}\bigr)\bigl(s+\frac{29}{12}\bigr)\bigl(s+\frac{5}{2}\bigr)\bigl(s+\frac{31}{12}\bigr)\bigl(s+\frac{11}{4}\bigr)\\
&&\bigl(s+\frac{35}{12}\bigr)\bigl(s+\frac{37}{12}\bigr), \\
b_{\a,x_3}(s)&=&\bigl(s+\frac{23}{12}\bigr)\bigl(s+2\bigr)\bigl(s+\frac{25}{12}\bigr)\bigl(s+\frac{9}{4}\bigr)\bigl(s+\frac{29}{12}\bigr)\bigl(s+\frac{5}{2}\bigr)\bigl(s+\frac{31}{12}\bigr)\bigl(s+\frac{11}{4}\bigr)\\
&&\bigl(s+\frac{17}{6}\bigr)\bigl(s+\frac{19}{6}\bigr), 
\end{eqnarray*}
and 
\begin{eqnarray*}
b_{\a,x_1^2}(s)&=&\bigl(s+2\bigr)\bigl(s+\frac{25}{12}\bigr)\bigl(s+\frac{9}{4}\bigr)\bigl(s+\frac{29}{12}\bigr)\bigl(s+\frac{5}{2}\bigr)\bigl(s+\frac{31}{12}\bigr)\bigl(s+\frac{11}{4}\bigr)\bigl(s+\frac{17}{6}\bigr)\\
&&\bigl(s+\frac{35}{12}\bigr)\bigl(s+\frac{19}{6}\bigr), \\
b_{\a,x_2^2}(s)&=&\bigl(s+2\bigr)\bigl(s+\frac{9}{4}\bigr)\bigl(s+\frac{29}{12}\bigr)\bigl(s+\frac{5}{2}\bigr)\bigl(s+\frac{31}{12}\bigr)\bigl(s+\frac{11}{4}\bigr)\bigl(s+\frac{17}{6}\bigr)\bigl(s+\frac{35}{12}\bigr)\\
&&\bigl(s+\frac{37}{12}\bigr)\bigl(s+\frac{19}{6}\bigr).  
\end{eqnarray*}
The ideal $J_{\f}(1)$ has a primary decomposition $\cap_{i=1}^{10}\q_i$ where  
\begin{eqnarray*}
\q_1&=&\langle 12s+17,x_1,x_2,x_3 \rangle,~ 
\q_2=\langle 2s+3,x_1,x_2,x_3 \rangle, \\
\q_3&=&\langle 12s+19,x_1,x_2,x_3 \rangle,~ 
\q_4=\langle 4s+7,x_1^2,x_2,x_3 \rangle,\\
\q_5&=&\langle 6s+11,x_1,x_2^2,x_3 \rangle, ~ 
\q_6=\langle 12s+23,x_1^2,x_1x_3,x_2,x_3^2 \rangle,\\
\q_7&=&\langle s+2,x_2^2-x_3x_1,x_1^3-x_3^2 \rangle, ~
\q_8=\langle 12s+25,x_1^3,x_1x_3,x_2,x_3^2 \rangle,\\
\q_9&=&\langle 6s+13,x_1^2,x_2^2,x_3 \rangle,~ 
\q_{10}=\langle 4s+9,x_1^3,x_1^2x_3,x_1x_2,x_1x_3-x_2^2,x_2x_3,x_3^2 \rangle,  
\end{eqnarray*}
hence the multiplier ideals are
\begin{equation*}
\J(\a^c)=
\left\{
	\begin{array}{lll}
	\C[x_1,x_2,x_3] & \mbox{ $0\le c<\frac{17}{12}$,} \\[1mm]
	\langle x_1, x_2,x_3 \rangle & \mbox{ $\frac{17}{12}\le c<\frac{7}{4}$,} \\[1mm]
	\langle x_1^2,x_2,x_3 \rangle & \mbox{ $\frac{7}{4}\le c<\frac{11}{6}$,}\\[1mm]
	\langle x_1^2,x_1x_2, x_2^2,x_3 \rangle & \mbox{ $\frac{11}{6}\le c<\frac{23}{12}$,}\\[1mm]
	\langle x_1,x_2,x_3 \rangle^2 & \mbox{ $\frac{23}{12}\le c<2$,}
	\end{array}
\right. 
\end{equation*}
and $\J(\a^c)=\a\J(\a^{c-1})$ for $c\ge 2$ as $\mu(\a)=2$. 
Note that $\frac{17}{12}$, $\frac{7}{4}$, $\frac{11}{6}$, $\frac{23}{12}$, $2$ are all jumping coefficients in $(0,2]$, 
and $\frac{3}{2}$, $\frac{19}{12}$, $\frac{25}{12}$, $\frac{13}{6}$, and $\frac{9}{4}$ are roots of $b_\a(-s)$ not coming from the jumping coefficients. 
\end{Example}
\begin{Example}
Let $\a=\langle x_1^3-x_2^2, x_2^3-x_3^2\rangle\subset \C[x_1,x_2,x_3]$ be the defining ideal of the space monomial curve $\Spec\C[T^4,T^6,T^9]\subset \C^3$. Then 
\begin{eqnarray*}
b_\a(s)&=&\bigl(s+\frac{4}{3}\bigr)\bigl(s+\frac{25}{18}\bigr)\bigl(s+\frac{29}{18}\bigr)\bigl(s+\frac{5}{3}\bigr)\bigl(s+\frac{31}{18}\bigr)\bigl(s+\frac{11}{6}\bigr)\bigl(s+\frac{35}{18}\bigr)\bigl(s+2\bigr)\\
&& \bigl(s+\frac{37}{18}\bigr)\bigl(s+\frac{13}{6}\bigr)\bigl(s+\frac{41}{18}\bigr), \\
b_{\a,x_1}(s)&=&\bigl(s+\frac{29}{18}\bigr)\bigl(s+\frac{5}{3}\bigr)\bigl(s+\frac{11}{6}\bigr)\bigl(s+\frac{35}{18}\bigr)\bigl(s+2\bigr)\bigl(s+\frac{37}{18}\bigr)\bigl(s+\frac{13}{6}\bigr)\bigl(s+\frac{41}{18}\bigr)\\
&& \bigl(s+\frac{7}{3}\bigr)\bigl(s+\frac{43}{18}\bigr)\bigl(s+\frac{49}{18}\bigr)\bigl(s+\frac{17}{6}\bigr), \\
b_{\a,x_2}(s)&=&\bigl(s+\frac{31}{18}\bigr)\bigl(s+\frac{35}{18}\bigr)\bigl(s+2\bigr)\bigl(s+\frac{37}{18}\bigr)\bigl(s+\frac{13}{6}\bigr)\bigl(s+\frac{41}{18}\bigr)\bigl(s+\frac{7}{3}\bigr)\bigl(s+\frac{43}{18}\bigr)\\
&& \bigl(s+\frac{47}{18}\bigr)\bigl(s+\frac{8}{3}\bigr)\bigl(s+\frac{17}{6}\bigr)\bigl(s+\frac{19}{6}\bigr),
\end{eqnarray*}
and the multiplier ideals are 
\begin{equation*}
\J(\a^c)=
\left\{
	\begin{array}{lll}
	\C[x_1,x_2,x_3] & \mbox{ $0\le c<\frac{4}{3}$,} \\[1mm]
	\langle x_1, x_2,x_3 \rangle & \mbox{ $\frac{4}{3}\le c<\frac{29}{18}$,} \\[1mm]
	\langle x_1^2,x_2,x_3 \rangle & \mbox{ $\frac{29}{18}\le c<\frac{31}{18}$,}\\[1mm]
	\langle x_1^2,x_1x_2, x_2^2,x_3 \rangle & \mbox{ $\frac{31}{18}\le c<\frac{11}{6}$,}\\[1mm]
	\langle x_1^3,x_1x_2,x_2^2,x_1x_3,x_2x_3,x_3^2 \rangle & \mbox{ $\frac{11}{6}\le c<\frac{35}{18}$,}\\[1mm]
	\langle x_1^3,x_1^2x_2,x_2^2,x_1x_3,x_2x_3,x_3^2 \rangle & \mbox{ $\frac{35}{18}\le c<2$,}
	\end{array}
\right. 
\end{equation*}
and $\J(\a^c)=\a\J(\a^{c-1})$ for $c\ge 2$. 
\end{Example}
\begin{Example}
Let $\a=\langle x_1^3-x_2^2, x_3^2-x_1^2x_2\rangle\subset \C[x_1,x_2,x_3]$ be the defining ideal of the space monomial curve $\Spec\C[T^4,T^6,T^7]\subset \C^3$. Then 
\begin{eqnarray*}
b_\a(s)&=&\bigl(s+\frac{4}{3}\bigr)\bigl(s+\frac{19}{14}\bigr)\bigl(s+\frac{23}{14}\bigr)\bigl(s+\frac{5}{3}\bigr)\bigl(s+\frac{25}{14}\bigr)\bigl(s+\frac{11}{6}\bigr)\bigl(s+\frac{27}{14}\bigr)\bigl(s+2\bigr)\\
&& \bigl(s+\frac{29}{14}\bigr)\bigl(s+\frac{13}{6}\bigr)\bigl(s+\frac{31}{14}\bigr), 
\end{eqnarray*}
and the multiplier ideals are 
\begin{equation*}
\J(\a^c)=
\left\{
	\begin{array}{lll}
	\C[x_1,x_2,x_3] & \mbox{ $0\le c<\frac{4}{3}$,} \\[1mm]
	\langle x_1, x_2,x_3 \rangle & \mbox{ $\frac{4}{3}\le c<\frac{23}{14}$,} \\[1mm]
	\langle x_1^2,x_2,x_3 \rangle & \mbox{ $\frac{23}{14}\le c<\frac{25}{14}$,}\\[1mm]
	\langle x_1^2,x_1x_2, x_2^2,x_3 \rangle & \mbox{ $\frac{25}{14}\le c<\frac{11}{6}$,}\\[1mm]
	\langle x_1,x_2,x_3 \rangle^2 & \mbox{ $\frac{11}{6}\le c<\frac{27}{14}$,}\\[1mm]
	\langle x_1^3,x_1x_2,x_1x_3,x_2^2,x_2x_3,x_3^2 \rangle & \mbox{ $\frac{27}{14}\le c<2$,}
	\end{array}
\right. 
\end{equation*}
and $\J(\a^c)=\a\J(\a^{c-1})$ for $c\ge 2$. 
\end{Example}
\begin{Example}
Let $\a=\langle x_1^4-x_2^3,x_3^2-x_1x_2^2\rangle\subset \C[x_1,x_2,x_3]$ be the defining ideal of the space monomial curve $\Spec\C[T^6,T^8,T^{11}]\subset \C^3$. Then 
\begin{eqnarray*}
b_\a(s)\!\!&=&\!\!
\bigl(s+\frac{9}{8}\bigr)^2\bigl(s+\frac{29}{24}\bigr)\bigl(s+\frac{31}{24}\bigr)\bigl(s+\frac{11}{8}\bigr)^2\bigl(s+\frac{35}{24}\bigr)\bigl(s+\frac{3}{2}\bigr)\bigl(s+\frac{37}{24}\bigr)\bigl(s+\frac{19}{12}\bigr)\\
&&\bigl(s+\frac{13}{8}\bigr)^2\bigl(s+\frac{41}{24}\bigr)\bigl(s+\frac{43}{24}\bigr)\bigl(s+\frac{11}{6}\bigr)\bigl(s+\frac{15}{8}\bigr)^2\bigl(s+\frac{23}{12}\bigr)\bigl(s+\frac{47}{24}\bigr)\bigl(s+2\bigr)\\
&&\bigl(s+\frac{49}{24}\bigr)\bigl(s+\frac{25}{12}\bigr)\bigl(s+\frac{13}{6}\bigr)\bigl(s+\frac{29}{12}\bigr), 
\end{eqnarray*}
and the multiplier ideals are 
\begin{equation*}
\J(\a^c)=
\left\{
	\begin{array}{lll}
	\C[x_1,x_2,x_3] & \mbox{ $0\le c<\frac{9}{8}$,} \\[1mm]
	\langle x_1, x_2,x_3 \rangle & \mbox{ $\frac{9}{8}\le c<\frac{11}{8}$,} \\[1mm]
	\langle x_1^2,x_2,x_3 \rangle & \mbox{ $\frac{11}{8}\le c<\frac{35}{24}$,}\\[1mm]
	\langle x_1^2,x_1x_2,x_2^2,x_3 \rangle & \mbox{ $\frac{35}{24}\le c<\frac{19}{12}$,}\\[1mm]
	\langle x_1^2,x_1x_2,x_2^2,x_1x_3,x_2x_3,x_3^2 \rangle & \mbox{ $\frac{19}{12}\le c<\frac{13}{8}$,}\\[1mm]
	\langle x_1^3,x_1x_2,x_2^2,x_1x_3,x_2x_3,x_3^2 \rangle & \mbox{ $\frac{13}{8}\le c<\frac{41}{24}$,}\\[1mm]
	\langle x_1^3,x_1^2x_2,x_2^2,x_1x_3,x_2x_3,x_3^2 \rangle & \mbox{ $\frac{41}{24}\le c<\frac{43}{24}$,}\\[1mm]
	\langle x_1^3,x_1^2x_2,x_1x_2^2,x_2^3,x_1x_3,x_2x_3,x_3^2 \rangle & \mbox{ $\frac{43}{24}\le c<\frac{11}{6}$,} \\[1mm]
	\langle x_1^3,x_1^2x_2,x_1x_2^2,x_2^3,x_1^2x_3,x_2x_3,x_3^2 \rangle & \mbox{ $\frac{11}{6}\le c<\frac{15}{8}$,}\\[1mm]
	\langle x_1^4,x_1^2x_2,x_1x_2^2,x_2^3,x_1^2x_3,x_2x_3,x_3^2 \rangle & \mbox{ $\frac{15}{8}\le c<\frac{23}{12}$,}\\[1mm]
	\langle x_1^4,x_1^2x_2,x_1x_2^2,x_2^3,x_1^2x_3,x_1x_2x_3,x_2^2x_3,x_3^2 \rangle & \mbox{ $\frac{23}{12}\le c<\frac{47}{24}$,}\\[1mm]
	\langle x_1^4,x_1^3x_2,x_1x_2^2,x_2^3,x_1^2x_3,x_1x_2x_3,x_2^2x_3,x_3^2 \rangle & \mbox{ $\frac{47}{24}\le c<2$,}
	\end{array}
\right. 
\end{equation*}
and $\J(\a^c)=\a\J(\a^{c-1})$ for $c\ge 2$. 
\end{Example}
\begin{Example}
Let 
$\a=\langle x_1^2-x_2x_3, x_2^2-x_1x_3, x_3^2-x_1x_2\rangle=I_2\left(
	\begin{array}{ccc}
	x_1 & x_2 & x_3 \\
	x_3 & x_1 & x_2 
	\end{array}
\right)\subset \C[x_1,x_2,x_3]$. 
Then 
\begin{eqnarray*}
b_\a(s)&=&\bigl(s+\frac{3}{2}\bigr)\bigl(s+ 2 \bigr)^2, \\
b^{(2)}_{\a}(s)&=&\bigl(s+\frac{3}{2}\bigr)\bigl(s+ 2 \bigr)^2(s+\frac{5}{2}\bigr)\bigl(s+ 3 \bigr)^2, \\
b_{\a,x_i}(s)&=&\bigl(s+ 2 \bigr)^2\bigl(s+\frac{5}{2}\bigr), 
\end{eqnarray*} 
and the multiplier ideals are
\begin{equation*}
J(f^c)=
\left\{
	\begin{array}{lll}
	\C[x_1,x_2,x_3] & \mbox{ $0\le c<\frac{3}{2}$,} \\[1mm]
	\langle x_1, x_2,x_3 \rangle & \mbox{ $\frac{3}{2}\le c<2$,}\\[1mm]
	\a & \mbox{$2 \le c<\frac{5}{2}$,}\\[1mm]
	\langle x_1, x_2,x_3 \rangle\a & \mbox{$\frac{5}{2}\le c<3$,}
	\end{array}
\right.
\end{equation*}
and $\J(\a^c)=\a\J(\a^{c-1})$ for $c\ge 3$. 
\end{Example}
\begin{Example}
Let $\a=\langle x_1^3-x_2x_3, x_2^2-x_1x_3, x_3^2-x_1^2x_2\rangle\subset \C[x_1,x_2,x_3]$ be the defining ideal of the space monomial curve $\Spec\C[T^3,T^4,T^5]\subset \C^3$. 
Then 
\[
b_\a(s)=\bigl(s+\frac{13}{9}\bigr)\bigl(s+\frac{3}{2}\bigr)\bigl(s+\frac{14}{9}\bigr)\bigl(s+\frac{16}{9}\bigr)\bigl(s+\frac{17}{9}\bigr)\bigl(s+2\bigr)^2\bigl(s+\frac{19}{9}\bigr)\bigl(s+\frac{20}{9}\bigr), 
\]
and 
\begin{equation*}
\J(\a^c)=
\left\{
	\begin{array}{lll}
	\C[x_1,x_2,x_3] & \mbox{ $0\le c<\frac{13}{9}$,} \\[1mm]
	\langle x_1,x_2,x_3 \rangle & \mbox{ $\frac{13}{9}\le c<\frac{16}{9}$,} \\[1mm]
	\langle x_1^2,x_2,x_3 \rangle & \mbox{ $\frac{16}{9}\le c<\frac{17}{9}$,}\\[1mm]
	\langle x_1^2,x_1x_2, x_2^2,x_3 \rangle & \mbox{ $\frac{17}{9}\le c<2$.}\\[1mm]
	\a & \mbox{ $2 \le c<\frac{22}{9}$.}
	\end{array}
\right. 
\end{equation*}
We must compute the ideal $J_{\f}(2)$ in Theorem \ref{algorithm for computing multiplier ideals 1} to obtain multiplier ideals $\J(\a^c)$ for $c\ge \frac{22}{9}$, 
and to determine whether $\frac{22}{9}$ is a jumping coefficient or not. 
It is, however, a really hard computation.  
\end{Example}
\begin{Example}
Let $f=x^3z^3+y^3z^2+y^2\in \C[x,y,z]$ and $\a=\langle f \rangle$. 
Then $b_f(s)=\bigl(s+\frac{5}{6}\bigr)^2(s+1)\bigl(s+\frac{7}{6}\bigr)^2\bigl(s+\frac{3}{2}\bigr)$, and 
$b_{\a,x}(s)=\bigl(s+\frac{5}{6}\bigr)(s+1)\bigl(s+\frac{7}{6}\bigr)^2\bigl(s+\frac{3}{2}\bigr)\bigl(s+\frac{11}{6}\bigr)$. 
The ideal $J_{f}(1)=I_{(f;1),2}$ has a primary decomposition $\cap_{i=1}^{8}\q_i$ where 
\begin{eqnarray*}
\q_1&=&\langle f,s+1 \rangle, \\
\q_2&=&\langle x,y,6s+5 \rangle, ~~~
\q_3 = \langle x^2,y,6s+7\rangle, \\
\q_4&=&\langle y,z,6s+5 \rangle, ~~~
\q_5 = \langle y,z^2,6s+7 \rangle, \\
\q_6&=&\langle x^3,y,z^3,xz,(6s+5)z,(6s+5)x,(6s+5)^2 \rangle, ~~~
\q_7 = \langle x^3,y,z^3,2s+3 \rangle, \\
\q_8&=&\langle x^3,y,z^3,x^2z^2,(6s+7)z^2,(6s+7)x^2,(6s+7)^2 \rangle. 
\end{eqnarray*}
The ideal $I_{(f;x),2}$ has primary decomposition 
$\cap_{i=1}^{7}\q'_i$ where 
\begin{eqnarray*}
\q'_1&=&\langle f,s+1 \rangle, \\
\q'_2&=&\langle x,y,6s+7 \rangle, ~~~
\q'_3 = \langle x^3,y,6s+11 \rangle, \\
\q'_4&=&\langle y,z,6s+5 \rangle, ~~~
\q'_5 = \langle y,z^2,6s+7 \rangle, \\
\q'_6&=&\langle x^3,y,z^3,xz^2,(6s+7)^2,(6s+7)x,(6s+7)z^2 \rangle, ~~~
\q'_7 = \langle x^2,y,z^3,2s+3 \rangle. 
\end{eqnarray*}
Let $\p_1=\langle x,y\rangle$ and $\p_2=\langle y,z\rangle$ be prime ideals of $\C[x,y,z]$, and let $\m_0=\langle x,y,z \rangle$ be the maximal ideal at the origin. 
Then we obtain local generalized Bernstein-Sato polynomials 
\begin{eqnarray*}
b^\a_f(s)=s+1, ~
b^{\p_1}_f(s)=b^{\p_2}_f(s)=\bigl(s+\frac{5}{6}\bigr)(s+1)\bigl(s+\frac{7}{6}\bigr), ~
b^{\m_0}_f(s)=b_f(s), 
\end{eqnarray*}
$b^\m_f(s)=\bigl(s+\frac{5}{6}\bigr)(s+1)\bigl(s+\frac{7}{6}\bigr)$ for $\m\in V(\p_1)\cup V(\p_2)\backslash \{\m_0\}$, 
$b^\p_f(s)=s+1$ for $\p \in V(\a)\backslash (V(\p_1)\cup V(\p_2))$, 
and $b^\p_f(s)=1$ for $\p \not\in V(\a)$, and $b^\a_{\a,x}(s)=s+1$, 
\begin{eqnarray*}
b^{\p_1}_{\a,x}(s)=(s+1)\bigl(s+\frac{7}{6}\bigr)\bigl(s+\frac{11}{6}\bigr), ~
b^{\p_2}_{\a,x}(s)=\bigl(s+\frac{5}{6}\bigr)(s+1)\bigl(s+\frac{7}{6}\bigr), 
\end{eqnarray*}
$b^{\m_0}_{\a,x}(s)=b_{\a,x}(s)$, $b^\m_{\a,x}(s)=b^{\p_i}_{\a,x}$ for $\m\in V(\p_i)\backslash \{\m_0\}$, 
$b^\p_{\a,x}(s)=s+1$ for $\p \in V(\a)\backslash (V(\p_1)\cup V(\p_2))$, and $b^\p_{\a,x}(s)=1$ for $\p \not\in V(\a)$. 
\end{Example}
\begin{Example}
Let $\a$ be a ideal of $\C[x,y,z]$ with a system of generators $\f=(x^3-y^2z,x^2+y^2+z^2-1)$. Then $b_\a(s)=\bigr(s+\frac{11}{6}\bigl)(s+2)\bigr(s+\frac{13}{6}\bigl)$. 
The ideal $J_{\f}(1)=I_{(\f;1),2}$ is 
\[
\langle \a, b_\a(s),(s+2)y,(s+2)(6s+13)x,(s+2)(z+1)(z-1)\rangle, 
\]
and its primary decomposition is $\q_1\cap\q_2\cap\q_3\cap\q_4\cap\q_5$ where 
\begin{eqnarray*}
\q_1&=&\langle \a, s+2 \rangle,\\
\q_2&=&\langle x,y,z-1,6s+11 \rangle, ~~~
\q_3 = \langle x,y,z-1,6s+13 \rangle,\\
\q_4&=&\langle x,y,z+1,6s+11 \rangle, ~~~
\q_5 = \langle x,y,z+1,6s+13 \rangle. 
\end{eqnarray*}
Thus 
\begin{equation*}
\J(\a^c)=
\left\{
	\begin{array}{lll}
	\C[x,y,z] & \mbox{ $0\le c<\frac{11}{6}$,} \\[1mm]
	\langle x,y,z^2-1 \rangle & \mbox{ $\frac{11}{6}\le c<2$.}
	\end{array}
\right. 
\end{equation*}
Let $\m_1=\langle x,y,z-1 \rangle$ and $\m_2=\langle x,y,z+1\rangle$ be maximal ideals of $\C[x,y,z]$. 
Then we obtain local generalized Bernstein-Sato polynomials 
\[
b^\a_\a(s)=s+2, ~~
b^{\m_1}_\a(s)=b^{\m_2}_\a(s)=\bigr(s+\frac{11}{6}\bigl)(s+2)\bigr(s+\frac{13}{6}\bigl), 
\]
$b^\m_{\a}(s)=s+2$ for $\m\in V(\a)\backslash \{\m_1,\m_2\}$, and $b^\p_{\a}(s)=1$ for $\p \not\in V(\a)$. 
\end{Example}

\end{document}